\documentclass[12pt]{article}

\usepackage{amsthm,amsmath,amssymb,accents}

\date{}

\iftrue \makeatletter
\@ifundefined{every@math@size}{}{%
\addto@hook\every@math@size{\dch@scr@hook}
\def\dch@scr@adjust{\@ifundefined{dch@sizet\f@size}%
  {\expandafter\dch@set@script\csname dch@sizet\f@size\endcsname}%
  {\csname dch@sizet\f@size\endcsname}}
\def\dch@set@script#1{\begingroup %
  \frozen@everymath{}%
  \let#1\@empty \let\dch@do@one\relax
  \dch@set@one\scriptscriptstyle\scriptscriptfont#1\ssf@size
  \dch@set@one\scriptstyle\scriptfont#1\sf@size
  \dch@set@one\textstyle\textfont#1\f@size
  \endgroup #1} %
\def\dch@set@one#1#2#3#4{%
  \@ifundefined{dch@size#4}%
   {\expandafter\xdef\csname dch@size#4\endcsname{%
      \fontdimen13\the#2\tw@=\the\fontdimen13#2\tw@
      \fontdimen14\the#2\tw@=\the\fontdimen14#2\tw@
      \fontdimen15\the#2\tw@=\the\fontdimen15#2\tw@
      \fontdimen16\the#2\tw@=\the\fontdimen16#2\tw@
      \fontdimen17\the#2\tw@=\the\fontdimen17#2\tw@}%
  }{\csname dch@size#4\endcsname}%
  \setbox\z@\hbox{$#1H_2$}\@tempdima\dp\z@
  \setbox\z@\hbox{$#1H_2^{+\vrule \@height 1em}$}%
   \ifdim\@tempdima<\dp\z@
    \advance\@tempdima\dp\z@ \divide\@tempdima\tw@ %
    \@tempdimb\dp\z@ \advance\@tempdimb-\@tempdima %
    \advance\@tempdimb\ht\z@ \advance\@tempdimb-1em %
    \xdef#3{#3\dch@do@one#2{\the\@tempdimb}{\the\@tempdima}}%
  \fi}
\def\dch@do@one#1#2#3{%
  \fontdimen13#1\tw@#2\relax
  \fontdimen14#1\tw@\fontdimen13#1\tw@ \fontdimen15#1\tw@\fontdimen13#1\tw@
  \fontdimen\sixt@@n#1\tw@#3\fontdimen17#1\tw@\fontdimen\sixt@@n#1\tw@}%
\let\dch@scr@hook\dch@scr@adjust
\ifx\glb@currsize\f@size \dch@scr@adjust \fi}%
\makeatother \fi
\normalsize

\newtheorem{thm}{Theorem}[section]
\newtheorem{cor}{Corollary}[section]
\newtheorem{lemma}{Lemma}[section]

\let\nfrac=\tfrac  % but use \tfrac
\def\dfrac#1#2{\lower0.15ex\hbox{\large$\frac{#1}{#2}$}}
\def\({\bigl(}
\def\){\bigr)}

\def\sumppd{\mathop{\sum\nolimits'}\limits}
\def\sumpp{\sum'}
\def\sumjkmn{\sum_{j=1}^m\sum_{k=1}^n}
\def\jsum{\sum_{j=1}^{m-1}}
\def\ksum{\sum_{k=1}^{n-1}}
\def\Amin{A_{\mathrm{min}}}

\def\b{{\mathord{\bullet}}}  %% full sum, j=1..m or k=1..n
\def\c{{\mathord{\ast}}}       %% sum j=1..m-1 or k=1..n-1

\def\BB{G}               % now old B = new G is a graph
\def\XX{H}  % redefine this now $X$ is a graph
\def\h{h}  %  h_{jk} gives edges of H.  Used to be x_{jk}

\def\II{\boldsymbol{I}}
\def\VV{\boldsymbol{V}}
\def\UU{\boldsymbol{U}}
\def\YY{\boldsymbol{Y}}
\def\MM{\boldsymbol{M}}
\def\barXX{{\mkern 5mu\overline{\mkern-5mu\XX\mkern-3mu}\mkern 3mu}}
\def\Vd{\VV_{\!\!\rm d}^{}}
\def\Vdinv{\VV_{\!\!\rm d}^{-1}}
\def\Vnd{\VV_{\!\!\rm nd}^{}}

\def\B{{\cal B}}
\def\D{{\cal D}}
\def\I{{\cal I}}

\def\M{{\cal M}}

\def\Q{{\cal Q}}
\def\R{{\cal R}}
\def\S{{\cal S}}

\def\W{{\cal W}}
\def\X{{\cal X}}
\def\Y{{\cal Y}}

\let\eps=\varepsilon
\let\t=\theta
\let\p=\phi

\def\varnu{{\dot\nu}}

\def\svec{\boldsymbol{s}}
\def\tvec{\boldsymbol{t}}
\def\vvec{\boldsymbol{v}}
\def\xvec{\boldsymbol{x}}
\def\yvec{\boldsymbol{y}}
\def\zvec{\boldsymbol{z}}

\def\betavec{\boldsymbol{\betavec}}
\def\thetavec{\boldsymbol{\theta}}
\def\phivec{\boldsymbol{\phi}}
\def\sigmavec{\boldsymbol{\sigma}}
\def\tauvec{\boldsymbol{\tau}}
\def\xivec{\boldsymbol{\xi}}
\def\zetavec{\boldsymbol{\zeta}}
\def\barX{{\mkern 5mu\overline{\mkern-5mu X\mkern-3mu}}}
\def\barY{{\mkern 2mu\overline{\mkern-2mu Y\mkern-2mu}}}

\def\avtheta{\bar{\thetavec}}
\def\avphi{\bar{\phivec}}
\def\avv{\bar{\vvec}}

\def\Chat{\widehat C}
\def\thetahatvec{\hat{\thetavec}}
\def\phihatvec{\hat{\phivec}}
\def\that{\hat\theta}  
\def\phat{\hat\phi}    
\def\varthetahatvec{\hat{\boldsymbol{\vartheta}}}
\def\varphihatvec{\hat{\boldsymbol{\varphi}}}

\def\Xtwo{R_{0,2}}
\def\Ytwo{C_{0,2}}
\def\Xthree{R_{0,3}}
\def\Ythree{C_{0,3}}
\def\pST{Y}
\def\XS{R_{1,1}}
\def\YT{C_{1,1}}
\def\XStwo{R_{2,1}}
\def\YTtwo{C_{2,1}}
\def\VSone{R_{1,2}}
\def\VTone{C_{1,2}}
\def\Pfour{Y_{1,1}}
\def\XtwoS{Y_{0,1}}
\def\YtwoT{Y_{1,0}}
\def\pSk{K}
\def\pTj{J}

\def\dsx{\delta}
\def\dty{\eta}

\def\OO{\widetilde{O}}

\def\Prob{\operatorname{Prob}}

\def\expect{\operatorname{\mathbb E}}
\def\hit{\operatorname{hit}}
\def\miss{\operatorname{miss}}
\def\aut{\operatorname{aut}}
\def\abs#1{\lvert#1\rvert} \let\card=\abs
\DeclareMathOperator{\tr}{tr}
\def\nicebreak{\vskip 0pt plus 50pt\penalty-300\vskip 0pt plus -50pt }
\def\ssqrt#1{\sqrt{\vrule width0pt height1.3ex depth0.4ex
  \smash[b]{#1}}} % a bit smashed

\begin{document}

\title{Random dense bipartite 
graphs and directed\\ graphs with specified degrees\thanks{
This is a revised expanded version of \textit{``Asymptotic enumeration
of dense 0-1 matrices with specified line sums and forbidden positions''},
submitted to arXiv on 22 Jan 2007.  A concise version will appear in
Random Structures and Algorithms.}
}
\author{
Catherine~Greenhill\\
\small School of Mathematics and Statistics\\[-0.8ex]
\small University of New South Wales\\[-0.8ex]
\small Sydney, Australia 2052\\[-0.3ex]
\small\texttt{csg@unsw.edu.au}
\and
\vrule width0pt height3ex
Brendan~D.~McKay\vrule width0pt height2ex\thanks
 {Research supported by the Australian Research Council.}\\
\small Department of Computer Science\\[-0.8ex]
\small Australian National University\\[-0.8ex]
\small Canberra ACT 0200, Australia\\[-0.3ex]
\small\texttt{bdm@cs.anu.edu.au}
}

\maketitle

\begin{abstract}
Let $\svec$ and $\tvec$ be vectors of positive integers with the same sum.
We study the uniform distribution on the space of simple bipartite graphs
with degree sequence~$\svec$ in one part and $\tvec$ in the other;
equivalently, binary matrices with row sums~$\svec$ and column
sums~$\tvec$.
In particular, we find precise formulae for the probabilities
that a given bipartite graph is edge-disjoint
from, a subgraph of, or an induced subgraph of a random graph 
in the class.
We also give similar formulae for the uniform distribution on
the set of simple directed graphs with out-degrees $\svec$ and
in-degrees~$\tvec$.
In each case, the graphs or digraphs are required to be sufficiently
dense, with the degrees varying within certain limits, and the
subgraphs are required to be sufficiently sparse. Previous results
were restricted to spaces of sparse graphs. Our theorems are
based on an enumeration of bipartite graphs avoiding a given set
of edges, proved by multidimensional complex integration. As a
sample application, we determine the expected permanent of a random
binary matrix with row sums $\svec$ and column sums $\tvec$.
\end{abstract}

\nicebreak
\section{Introduction}\label{s:intro}

Let $\svec=(s_1,\ldots, s_m)$ and
$\tvec = (t_1,\ldots, t_n)$ be vectors of positive integers
with  $\sum_{j=1}^{m}s_j = \sum_{k=1}^n t_k$.  Define
$\B(\svec,\tvec)$ to be the set of simple bipartite graphs
with vertices $\{ u_1,\ldots, u_m\}\cup \{ v_1,\ldots, v_n\}$,
such that vertex $u_j$ has degree $s_j$ for $j=1,\ldots, m$
and vertex $v_k$ has degree~$t_k$ for $k=1,\ldots, n$.
Equivalently, we may think of $\B(\svec,\tvec)$ as the set
of all $m\times n$ matrices over $\{ 0,1\}$ with $j$th row sum equal
to $s_j$ for $j=1,\ldots, m$ and $k$th column sum equal to~$t_k$ for
$k=1,\ldots, n$.

In addition, let $\XX$ be a fixed bipartite graph on the same vertex
set. In this paper we find precise formulae for the probabilities
that $\XX$ is edge-disjoint from $\BB\in \B(\svec,\tvec)$, that $\XX$
is a subgraph of~$\BB$, and that $\XX$ is an induced subgraph of~$\BB$.
These probabilities are defined for the uniform distribution on
$\B(\svec,\tvec)$. In general, whenever we refer to a random element
of a set, we always mean an element chosen uniformly at random.

These formulae are obtained when the graphs in $\B(\svec,\tvec)$
are sufficiently dense, the graph $\XX$ is sufficiently sparse 
and the entries of $\svec$ and $\tvec$ only vary within certain limits.
The exact conditions are stated in Section~\ref{s:graphs}.
The starting point of the calculations is an enumeration of
the set $\B(\svec,\tvec,\XX)$ of graphs in $\B(\svec,\tvec)$
which are edge-disjoint from~$\XX$; see Theorem~\ref{bigtheorem}.

In the case $m=n$, the $n\times n$ binary matrix associated with the
bipartite graph can also be interpreted as the adjacency matrix
of a digraph which has no multiple edges but may have loops. By
excluding the diagonal we obtain a parallel series of results
for simple digraphs (digraphs without multiple edges or loops).
These are presented in Section~\ref{s:digraphs}.

These subgraph probabilities enable the development of a theory
of random graphs and digraphs in these classes. As examples of
computations made possible by this theory, we calculate the expected
number of subgraphs isomorphic to a given regular subgraph. A
particular case of interest is the permanent of a random 0-1 matrix
with row sums~$\svec$ and column sums~$\tvec$.

\medskip

Now we briefly review the history of this problem.
All previous precise asymptotics were restricted to
sparse graphs.
Define $g=\max\{s_1,\ldots,s_m,\allowbreak t_1,\ldots,t_n\}$,
$x=\max\{x_1,\ldots,x_m,\allowbreak y_1,\ldots,y_n\}$
and $N=\sum_j s_j$.
Asymptotic estimates for bounded $g$ were found by
Bender~\cite{Bender} and Wormald~\cite{Wormald}.
This was extended by Bollob\'as and McKay~\cite{Bollobas}
to the case $g,x=O\(\min\{\log m,\log n\}^{1/3}\)$ and
by McKay~\cite{Silver} to the case $g^2+xg=o(N^{1/2})$.
Estimates which are sometimes more widely applicable
were given by McKay~\cite{McKay81}.
The best enumerative results for $\B(\svec,\tvec)$
in the sparse domain appear in~\cite{GMW,MW}.

Although results about sparse digraphs with specified in-degree
and out-degree sequences can be deduced from the above, we are
not aware of this having been done.
Some results using the pairings model have appeared~\cite{CFM}.
For digraphs in the dense regime, some related work includes
enumeration of tournaments by score sequence with possible
forbidden subgraph~\cite{RT, MWtournament, eulerian, GaoTourn},
Eulerian digraphs~\cite{RT, Jimmy92}, Eulerian oriented
graphs~\cite{RT, Jimmy97}, and digraphs with a given excess
sequence~\cite{Jimmy95}.

For the case of dense bipartite graphs with specified degrees,
an asymptotic formula for the case of empty~$\XX$
was given by Canfield and McKay~\cite{CM} for semiregular graphs
and by Canfield, Greenhill and McKay~\cite{CGM} for irregular graphs.
The latter study is the inspiration for the present one.
A similar study for graphs which are not necessarily bipartite
is in preparation~\cite{RANX}.

In related work using different methods, Barvinok~\cite{barvinok}
gives upper and lower bounds for $\card{\B(\svec,\tvec,\XX)}$
which hold very generally (from sparse to dense graphs) but which
can differ by a factor of $(mn)^{O(m+n)}$.   Barvinok's results
also give insight into the structure of a ``typical'' element of
$\B(\svec,\tvec,\XX)$, which he proves is close to a certain
``maximum entropy'' matrix.

\medskip

The paper is structured as follows.
The results for bipartite graphs are presented in 
Section~\ref{s:graphs} and the corresponding results for digraphs
can be found in Section~\ref{s:digraphs}. 
Then Section~\ref{s:proof} presents a proof of the fundamental
enumeration
result, Theorem~\ref{bigtheorem}, from which everything else follows.

Throughout the paper, the asymptotic notation $O(f(m,n))$ refers
to the passage of~$m$ and~$n$ to $\infty$.
We also use a modified
notation $\OO(f(m,n))$, which is to be taken as a shorthand for
$O\(f(m,n)n^{O(1)\eps}\)$, where the $O(1)$ factor is uniform
over~$\eps$ provided $\eps$ is small enough.

\nicebreak
\section{Subgraphs of random bipartite graphs}\label{s:graphs}

In this section we state our results for bipartite graphs.

The starting point of the investigation is the enumeration formula
given in the following theorem.
Define $m,n,\svec,\tvec$ as in the
Introduction and further define
\[ s = m^{-1}\sum_{j=1}^m s_j, \quad
   t = n^{-1}\sum_{k=1}^n t_k, \quad
  \lambda = s/n=t/m, \quad
  A=\dfrac12\lambda(1-\lambda). 
 \]
Note that $s$ is the average degree on one side of the vertex
bipartition, $t$ is the average degree on the other side, and
$\lambda$ is the edge density (the number of edges divided by $mn$).  

Let $\XX$ be a fixed bipartite graph on the same vertex set that
defines $\B(\svec,\tvec)$, namely $\{u_1,\ldots,u_m\}\cup
\{v_1,\ldots,v_n\}$.
For $j=1,\ldots, m$, $k=1,\ldots, n$, let $x_j$ and $y_k$ be the
degrees of vertices $u_j$ and $v_k$ of $\XX$, respectively,
and further define
\begin{align*}
 && \dsx_j &= s_j-s+\lambda x_j, &
  \dty_k &= t_k-t+\lambda y_k. && 
\end{align*}
Also define
 \begin{align*}
  &&  X &= \sum_{j=1}^m x_j = \sum_{k=1}^n y_k, &
     \pST &= \sum_{jk\in\XX} \dsx_j \dty_k, &&\\
 && R &= \sum_{j=1}^m \,(s_j-s)^2, &
    C &= \sum_{k=1}^n \,(t_k-t)^2. &&
    \end{align*}
In the case of $Y$ and similar notation used in this section, the summation
is over all $j\in\{u_1,\ldots,u_m\}$
and $k\in\{v_1,\ldots,v_n\}$ such that $u_jv_k$ is an edge of~$\XX$.

\begin{thm}\label{bigtheorem}
For some $\eps>0$, suppose that
$s_j-s$, $x_j$, $t_k-t$ and $y_k$ are
uniformly $O(n^{1/2+\eps})$ for $1\le j\le m$ and $1\le k\le n$,
and $X=O(n^{1+2\eps})$, for $m,n\to\infty$.
Let\/ $a,b>0$ be constants such that\/ $a+b<\tfrac{1}{2}$.
Suppose that $m$, $n\to \infty$ with
$n = o(m^{1+\eps})$, 
$m=o(n^{1+\eps})$ and
\[ \frac{(1-2\lambda)^2}{8A}\biggl(1 + \frac{5m}{6n}+
   \frac{5n}{6m}\biggr) \leq a\log n. \]
Then, provided $\eps>0$ is small enough, we have
\begin{align*}
 \card{\B(\svec,\tvec,\XX)} &= \binom{mn{-}X}{\lambda mn}^{\!\!-1}\,
\prod_{j=1}^m \binom{n{-}x_j}{s_j} \prod_{k=1}^n \binom{m{-}y_k}{t_k} \\
  &{\kern6mm}\times 
    \exp\biggl( -\frac12\Bigl(1-\frac{R}{2Amn}\Bigr)
    \Bigl(1-\frac{C}{2Amn}\Bigr) 
    - \frac{\pST}{2Amn} + O(n^{-b})\biggr).
\end{align*}
\end{thm}

\smallskip

The proof of Theorem~\ref{bigtheorem} will be presented in
Section~\ref{s:proof}.  As in the special case of empty $H$
proved in~\cite{CGM}, the formula for $\card{\B(\svec,\tvec,\XX)}$
has an intuitive interpretation. The first binomial and
the two products of binomials are,
respectively, the number of graphs with $\lambda mn$ edges
that avoid $H$, the number of such graphs with row sums~$\svec$,
and the number of such graphs with column sums~$\tvec$.
Therefore, the exponential factor measures the non-independence of
the events of having row sums~$\svec$ and having column sums~$\tvec$.
Another expression for the product of binomials in the theorem is
given below in equation~\eqref{rad16}.

\bigskip

We can now employ
Theorem~\ref{bigtheorem} to explore the uniform probability
space over $\B(\svec,\tvec)$.
First we need a little more notation.
For all nonnegative integers $h,\ell$ define
\begin{align*}
 && R_{h,\ell} &= \sum_{j=1}^m \dsx_j^h x_j^\ell, &
  C_{h,\ell} &= \sum_{k=1}^n \dty_k^h y_k^\ell. && 
\end{align*}
We will abbreviate $R_{h,0}=R_h$ and $C_{h,0}=C_h$.
Also note that $R_1=C_1=\lambda X$ and $R_{0,1}=C_{0,1}=X$.
Finally, let
\begin{align*}
 \Pfour = \sum_{jk\in\XX} x_jy_k,   \qquad
\XtwoS = \sum_{jk\in\XX} \dsx_jy_k,  \qquad
   \YtwoT = \sum_{jk\in\XX} x_j\dty_k. 
\end{align*}

\begin{thm}\label{avoidingXbip}
Under the conditions of Theorem~\ref{bigtheorem}, 
the following are true for a random graph $\BB\in\B(\svec,\tvec)$ 
provided $\eps>0$ is small enough:
\begin{enumerate}
\itemsep=0pt
\item 
the probability that $\BB$
is edge-disjoint from $\XX$ is\/ $(1-\lambda)^X \miss(m,n)$;
\item 
the probability that $\BB$ contains $\XX$ as a subgraph is\/
$\lambda^X \hit(m,n),$
\end{enumerate}
where
 \begin{align*}
   \miss(m,n) &= \exp\biggl(
     \frac{\lambda X}{2(1-\lambda)}\Bigl(\frac{1}{n}+\frac{1}{m}\Bigr)
     + \frac{\lambda X^2}{2(1-\lambda)mn}
     -\frac{1}{1-\lambda}\Bigl(\frac{\XS}{n}+\frac{\YT}{m}\Bigr) \\
    &{\kern14mm}
     -\frac{\pST}{\lambda(1-\lambda)mn}
     +\frac{\lambda}{2(1-\lambda)}
        \Bigl(\frac{\Xtwo}{n}+\frac{\Ytwo}{m}\Bigr)        
     +\frac{\lambda(1-2\lambda)}{6(1-\lambda)^2}
        \Bigl(\frac{\Xthree}{n^2}+\frac{\Ythree}{m^2}\Bigr) \\
    &{\kern14mm}
     -\frac{1-2\lambda}{2(1-\lambda)^2}
        \Bigl(\frac{\VSone}{n^2}+\frac{\VTone}{m^2}\Bigr)
     -\frac{1}{2(1-\lambda)^2}
        \Bigl(\frac{\XStwo}{n^2}+\frac{\YTtwo}{m^2}\Bigr)
     + O(n^{-b})\biggr)
 \end{align*}
 and
 \begin{align*}
   \hit(m,n) &= \exp\biggl(
     \frac{(1-\lambda) X}{2\lambda}
         \Bigl(\frac{1}{n}+\frac{1}{m}\Bigr)
     + \frac{(1-\lambda) X^2}{2\lambda mn}         
     + \frac{1}{\lambda}\Bigl(\frac{\XS}{n}+\frac{\YT}{m}\Bigr) \\
    &{\kern14mm}
     - \frac{1}{2\lambda^2}
        \Bigl(\frac{\XStwo}{n^2}+\frac{\YTtwo}{m^2}\Bigr)
     - \frac{1+\lambda}{2\lambda}
        \Bigl(\frac{\Xtwo}{n}+\frac{\Ytwo}{m}\Bigr)
     + \frac{1+2\lambda}{2\lambda^2}
        \Bigl(\frac{\VSone}{n^2}+\frac{\VTone}{m^2}\Bigr) \\
    &{\kern14mm}
     - \frac{(1+\lambda)(1+2\lambda)}{6\lambda^2}
        \Bigl(\frac{\Xthree}{n^2}+\frac{\Ythree}{m^2}\Bigr)
     - \frac{\pST-\XtwoS-\YtwoT+\Pfour}{\lambda(1-\lambda)mn}
     + O(n^{-b})\biggr).
 \end{align*}
\end{thm}

\begin{proof}  
 The first probability in the statement of Theorem~\ref{avoidingXbip} is
\[
    \frac{\card{\B(\svec,\tvec,\XX)}}{\card{\B(\svec,\tvec)}}
\]
which can be expanded using Theorem~\ref{bigtheorem}.
(One method is to apply~\eqref{rad16} below.)
The second probability can be derived in similar fashion,
or can be deduced from the first on noting that the
probability that $\BB$ includes~$\XX$ is the probability that
the complement of $\BB$ avoids~$\XX$.
\end{proof}
In the standard model of random bipartite graphs on $m+n$ vertices
with expected edge density~$\lambda$, each of the $mn$ possible
edges is present independently with probability~$\lambda$. 
The probability that a random bipartite graph taken from the
standard model is disjoint from or contains a given set of $X$ edges
is $(1-\lambda)^X$ or $\lambda^X$, respectively.
Therefore, the quantities $\miss(m,n)$ and
$\hit(m,n)$ given in
Theorem~\ref{avoidingXbip} can be interpreted as a measure of
how far these probabilities differ in $\B(\svec,\tvec)$ 
compared to the standard model. 
Suppose that in addition to the conditions of Theorem~\ref{avoidingXbip},
we also have
\begin{equation}
\begin{split}\label{biteme1}
  X\max_j\,\abs{s_j-s} + \lambda R_{0,2} &= o\((1-\lambda)n\), \\
  X\max_k\,\abs{t_k-t} + \lambda C_{0,2} &= o\((1-\lambda)m\). 
\end{split}
\end{equation}
Then $\miss(m,n)=1+o(1)$.
Similarly, if we have
\begin{equation}
\begin{split}\label{biteme2}
  X\max_j\,\abs{s_j-s} + (1-\lambda) R_{0,2} &= o(\lambda n), \\
  X\max_k\,\abs{t_k-t} + (1-\lambda) C_{0,2} &= o(\lambda m). 
\end{split}
\end{equation}
then $\hit(m,n)=1+o(1)$.
Requirements~\eqref{biteme1} and~\eqref{biteme2} are both met,
for example, if $X=O(n^{1/2-2\eps})$.
Another interesting case is when $s_j-s$, $x_j$, $t_k-t$
and $y_k$ are uniformly $O(n^\eps)$ and $X=O(n^{1-2\eps})$.

To assist with the application of Theorem~\ref{avoidingXbip}, we will
give the simplifications that result when the graphs in $\B(\svec,\tvec)$
are semiregular or when the graph $\XX$ is semiregular.

\begin{cor}\label{avoidingXflat}
In addition to the conditions of Theorem~\ref{avoidingXbip},
assume that $s_j=s$ and $t_k=t$ for all~$j,k$.
Then
 \begin{align*}
   \miss(m,n) &= \exp\biggl(  
     \frac{\lambda X}{2(1-\lambda)}\Bigl(\frac{1}{n}+\frac{1}{m}\Bigr)
     + \frac{\lambda X^2}{2(1-\lambda)mn}
     -\frac{\lambda\Pfour}{(1-\lambda)mn} \\
    &{\kern12mm}
     -\frac{\lambda}{2(1-\lambda)}
        \Bigl(\frac{\Xtwo}{n}+\frac{\Ytwo}{m}\Bigr)        
     -\frac{\lambda(2-\lambda)}{6(1-\lambda)^2}
        \Bigl(\frac{\Xthree}{n^2}+\frac{\Ythree}{m^2}\Bigr) 
     + O(n^{-b})\biggr)
 \end{align*}
 and 
 \begin{align*}
   \hit(m,n)  &= \exp\biggl(  
     \frac{(1-\lambda) X}{2\lambda}\Bigl(\frac{1}{n}+\frac{1}{m}\Bigr)
     + \frac{(1-\lambda) X^2}{2\lambda mn}
     -\frac{(1-\lambda)\Pfour}{\lambda mn} \\
    &{\kern12mm}
     -\frac{1-\lambda}{2\lambda}
        \Bigl(\frac{\Xtwo}{n}+\frac{\Ytwo}{m}\Bigr)        
     -\frac{1-\lambda^2}{6\lambda^2}
        \Bigl(\frac{\Xthree}{n^2}+\frac{\Ythree}{m^2}\Bigr) 
     + O(n^{-b})\biggr).\quad\qedsymbol
 \end{align*}
\end{cor}

\begin{cor}\label{avoidingXreg}
In addition to the conditions of Theorem~\ref{avoidingXbip},
assume that $x_j=x$ and $y_k=y$ for all~$j,k$.
(Note that Theorem~\ref{avoidingXbip} requires
$x,y=O(n^{2\eps})$ in that case.)
Then 
\[
  \miss(m,n) = \exp\biggl(  
     - \frac{\lambda(xy-x-y)}{2(1-\lambda)}
     - \frac{yR+xC}{2(1-\lambda)^2mn}
     - \frac{\widehat Y}{\lambda(1-\lambda)mn}
     + O(n^{-b})\biggr)
 \]
and 
 \[ \hit(m,n) =  \exp\biggl(  
     - \frac{(1-\lambda)(xy-x-y)}{2\lambda}
     - \frac{yR+xC}{2\lambda^2mn}
     - \frac{\widehat Y}{\lambda(1-\lambda)mn}
     + O(n^{-b})\biggr),
 \]
where $\widehat Y=\sum_{jk\in\XX}(s_j-s)(t_k-t)$.\quad\qedsymbol
\end{cor}

\medskip

The next question we will address is the probability of $\XX$
appearing as an induced subgraph.
To be precise, suppose that $\XX$ has no edges outside
$\{u_1,\ldots,u_J\}\times\{v_1,\ldots,v_K\}$ and let
$\XX_{J,K}$ denote the subgraph of $\XX$ induced by
those vertices. 
We will only consider the situation when the graphs
in $\B(\svec,\tvec)$ are semiregular.
The corresponding result for irregular graphs can also
be obtained using the same approach.

The probability that $\XX_{J,K}$ is an induced subgraph of
$\BB\in\B(\svec,\tvec)$ is simpler to state in terms of some
new variables.
For $\ell=1,2,3$, define
\[
 \omega_\ell = \sum_{j=1}^J (x_j-\lambda K)^\ell, \qquad
 \omega'_\ell = \sum_{k=1}^K (y_k-\lambda J)^\ell.
\]
Note that $\omega_1=\omega'_1=X-\lambda JK$.

\begin{thm}\label{inducedXbip}
Adopt the assumptions of Theorem~\ref{bigtheorem}
with $s_j=s$ and $t_k=t$ for all~$j,k$, and
assume that $J,K=O(n^{1/2+\eps})$. 
Then the probability that a random graph in $\B(\svec,\tvec)$
has $H_{J,K}$ as an induced subgraph is
\begin{align*}
    \lambda^X(1-\lambda)^{JK-X} \exp\biggl(&
    \Bigl(\frac{JK}{2}+\frac{(1-2\lambda)\omega_1}{4A}\Bigr)
          \Bigl(\frac{1}{m}+\frac{1}{n}\Bigr) - \frac{\omega_1^2}{4Amn}\\
    &{} -\frac{(n+K)\omega_2}{4An^2} 
     -\frac{(m+J)\omega'_2}{4Am^2}
     -\frac{1-2\lambda}{24A^2}
        \Bigl(\frac{\omega_3}{n^2}+\frac{\omega'_3}{m^2}\Bigr)
     + O(n^{-b})\biggr).
\end{align*}
\end{thm}

\begin{proof}
Let $\XX^\ast$ be the complete
bipartite graph on the parts $\{u_1,\ldots,u_J\}$ and $\{v_1,\ldots,v_K\}$.
Then the probability that a random graph in $\B(\svec,\tvec)$
has $\XX_{J,K}$ as an induced subgraph is
\[
    \frac{\card{\B(\svec-\xvec,\tvec-\yvec,\XX^\ast)}}
         {\card{\B(\svec,\tvec)}}.
\]
This ratio can be estimated using Theorem~\ref{bigtheorem} 
(or by combining Theorems~\ref{bigtheorem} and~\ref{avoidingXbip}).
\end{proof}

The argument of the exponential in Theorem~\ref{inducedXbip} is
$o(1)$ if $JK^2=o(An)$ and $J^2K=o(Am)$.
So, in those circumstances, the probabilities of induced subgraphs
asymptotically match the standard bipartite random graph model for
edge probability~$\lambda$.

\bigskip

A related question asks for the distribution of the number of
subgraphs of given type in a random graph in $\B(\svec,\tvec)$.
This deserves a serious study, which we will only just
initiate here.
A \textit{colour-preserving isomorphism\/} of two bipartite
graphs on $\{u_1,\ldots,u_m\}\cup\{v_1,\ldots,v_n\}$ is an
isomorphism that preserves the sets $\{u_1,\ldots,u_m\}$ and
$\{v_1,\ldots,v_n\}$.
Let $\I(\XX)$ be the set of all graphs isomorphic to $\XX$ by
a colour-preserving isomorphism.  We know that
\[
    \card{\I(\XX)} = \frac{m!\,n!}{\aut(\XX)},
\]
where $\aut(\XX)$ is the number of colour-preserving automorphisms
of~$\XX$.

When the graphs in $\B(\svec,\tvec)$ are semiregular,
the expected number of elements of $\I(\XX)$ that are
contained in or edge-disjoint from a random graph in
$\B(\svec,\tvec)$ is clearly just
$\card{\I(\XX)}$ times the probability given by
Theorem~\ref{avoidingXbip} and
Corollary~\ref{avoidingXflat}.

If this regularity condition does not hold, the calculation is
more complex.  Here we consider the case that the graph $\XX$ is
semiregular and leave the most general case for a future paper.

We will need the following averaging lemma.

\medskip
\begin{lemma}\label{averaging}
Let $\zvec^{(0)}=(z^{(0)}_1,z^{(0)}_2,\ldots,z^{(0)}_n)$ be a
vector in $[-1,1]^n$ such that $\sum_{j=1}^n z^{(0)}_j=0$.
Form $\zvec^{(1)}$, $\zvec^{(2)}, \ldots\,$ as follows:
for each $r\ge 0$,
if $z^{(r)}_i$ is the first of the smallest elements of
$\zvec^{(r)}$ and $z^{(r)}_\ell$ is the first of the 
largest elements of
$\zvec^{(r)}$, then $\zvec^{(r+1)}$ is the same as $\zvec^{(r)}$
except that $z^{(r+1)}_i$ and $z^{(r+1)}_\ell$ are both
equal to $(z^{(r)}_i+z^{(r)}_\ell)/2$.
Then $\zvec^{(n)}\in[-\tfrac12,\tfrac12]^n$.
\end{lemma}
\begin{proof}
If $\zvec^{(r)}\notin[-\tfrac12,\tfrac12]^n$, then the fact that
$\sum_{j=1}^n z^{(r)}_j=0$ implies that $z^{(r)}_i<0$
and $z^{(r)}_\ell>0$.  Therefore $\zvec^{(r+1)}$ has at least
one fewer element outside $[-\tfrac12,\tfrac12]$ than
$\zvec^{(r)}$ does.  The lemma follows.  (In fact,
$\zvec^{(\lfloor (2n-1)/3\rfloor)}\in[-\tfrac12,\tfrac12]^n$,
but this improvement is not necessary for our application.)
\end{proof}

\begin{thm}\label{unlabelled}
Suppose that the conditions of Theorem~\ref{bigtheorem} apply with
$x_j=x$ and $y_k=y$ for all $j,k$.
Then the following is true of a random graph $G$ in
$\B(\svec,\tvec)$:
\begin{enumerate}
\itemsep=0pt
\item the expected number of graphs in $\I(\XX)$ that are
subgraphs of $G$ is
\[
   \lambda^X\card{\I(\XX)}\exp\biggl(
      - \frac{(1-\lambda)(xy-x-y)}{2\lambda} 
       - \frac{yR+xC}{2\lambda^2mn}
      + O(n^{-b})\biggr);
\]
\item the expected number of graphs in $\I(\XX)$ that are
edge-disjoint from $G$ is
\[
   (1-\lambda)^X\card{\I(\XX)}\exp\biggl(
      - \frac{\lambda(xy-x-y)}{2(1-\lambda)} 
       - \frac{yR+xC}{2(1-\lambda)^2mn}
      + O(n^{-b})\biggr).
\]
\end{enumerate}
\end{thm}

\begin{proof}
Define $\zvec^{(0)},\zvec^{(1)},\ldots\,$ as in Lemma~\ref{averaging},
with $\zvec^{(0)}_j=s_j-s$ for $1\le j\le n$. 
For $r\ge 0$, define
\[
   Y^{(r)}(g,h) = \sum_{j=1}^m \,z^{(r)}_{j^g}\, T_{j,h},
   \text{ where ~}
   T_{j,h} = \sum_{k\,:\, u_jv_k\in E(\XX)} \!\!(t_{k^h}-t),
\]
and
\[
   F^{(r)} = \sum_{(g,h)\in S_m\times S_n}\!\!
      \exp\biggl( -\frac{\pST^{(r)}(g,h)}{2Amn} \biggr).
\]

For a permutation pair $(g,h)\in S_m\times S_n$, define
$\XX^{g,h}$ to be the isomorph of $\XX$ with edge set
$\{ u_{j^g}v_{k^h} \mid u_jv_k\in E(\XX)\}$.
As $(g,h)$ runs over $S_m\times S_n$, 
each isomorph of $\XX$ appears as $\XX^{g,h}$ 
exactly $\aut(\XX)$ times.  Therefore,
by Corollary~\ref{avoidingXreg}, the expectation required in
part~(i) of the theorem is
\[
   \frac{\lambda^X F^{(0)}}{\aut(\XX)}
   \exp\biggl(
      - \frac{(1-\lambda)(xy-x-y)}{2\lambda}
       - \frac{yR+xC}{2\lambda^2mn} + O(n^{-b})\biggr).
\]

For some $r\ge 0$,
suppose that $\zvec^{(r+1)}$ is formed from $\zvec^{(r)}$ by
averaging $z^{(r)}_i$ and $z^{(r)}_\ell$ as in Lemma~\ref{averaging}.
Then $\{ (i\ell)g \,|\, g\in S_n\}=S_n$, so
\begin{align*}
  F^{(r)} &= \dfrac12 \sum_{(g,h)} \biggl(
       \exp\Bigl( -\frac{Y^{(r)}(g,h)}{2Amn}  \Bigr)
       + \exp\Bigl( -\frac{Y^{(r)}((i\ell)g,h)}{2Amn}  
                       \Bigr)\biggr)\\
       &= \dfrac12 \sum_{(g,h)} 
            \exp\Bigl( -\frac{\sum_{j\notin\{i,\ell\}}z^{(r)}_j
                        T_{j,h}}{2Amn}  \Bigr)
         \biggl( \exp\Bigl(
     -\frac{z^{(r)}_i T_{i,h}+z^{(r)}_\ell T_{\ell,h}}{2Amn}
               \Bigr) \\
          &{\kern7.2cm}+ \exp\Bigl(
          -\frac{z^{(r)}_i T_{\ell,h}+z^{(r)}_\ell T_{i,h}}{2Amn}
              \Bigr)\biggr)
          \displaybreak[0]\\
      &= 
      \sum_{(g,h)}
           \exp\Bigl( -\frac{\sum_{j\notin\{i,\ell\}}z^{(r)}_jT_{j,h}}
            {2Amn}
           -\frac{z^{(r+1)}_i T_{i,h}+z^{(r+1)}_\ell T_{\ell,h}}
              {2Amn}
               + \OO(n^{-2})\Bigr) \\
      &= \sum_{(g,h)}
       \exp\Bigl( -\frac{Y^{(r+1)}(g,h)}{2Amn} 
                  + \OO(n^{-2}) \Bigr) \\
      &= F^{(r+1)}\exp\( \OO(n^{-2}) \).
\end{align*}
By Lemma~\ref{averaging} there is some $r_0=O(n\log n)$ such that
$\zvec^{(r_0)}\in[-n^{-1/2},n^{-1/2}]^n$.  By the definition
of $F^{(r_0)}$, we have
$F^{(r_0)}=m!\,n!\,\exp\( \OO(n^{-1}) \)$, so
$F^{(0)}=m!\,n!\,\exp\( \OO(n^{-1}) \)$ by induction.  Part~(i) of the
theorem follows.  Part~(ii) is proved in identical fashion.
\end{proof}

A simple example of Theorem~\ref{unlabelled} at work is 
the enumeration of perfect matchings in the case $m=n$.
Equivalently, this is the permanent of the corresponding
$n\times n$ binary matrix.
Most previous research has focussed on the case that the matrix
has constant row and column sums.  For $s=t=o(n^{1/3})$, the
asymptotic expectation and variance are known, while for
$s=t=n-O(n^{1-\epsilon})$, the asymptotic expectation is
known~\cite{Bollobas}.
In the intermediate range of densities covered by the current paper,
it appears that only bounds are known.  The van der Waerden lower
bound $n!\,\lambda^n$ (proved independently by Egorychev and
Falikman) was improved by Gurvits~\cite{Gurvits} to
$s!\,\((s-1)^{s-1}/s^{s-2}\)^{n-s}$.  The best upper bound
is $s!^{1/\lambda}\sim n!\,\lambda^{n+1/(2\lambda)}
(2\pi n)^{(1-\lambda)/(2\lambda)}$ conjectured by Minc and proved
by Bregman. 
See Timash\"ev~\cite{Timashev} for
references and discussion. 

Applying Theorem~\ref{unlabelled}(i) with $x=y=1$ gives the following.

\begin{thm}\label{permanent}
Suppose that $m=n$ and $\svec,\tvec,\lambda$ satisfy the
requirements of Theorem~\ref{bigtheorem}.
Then the expected
permanent of a random $n\times n$ matrix over $\{0,1\}$
with row sums $\svec$ and column sums $\tvec$ is
\[
   n!\,\lambda^n\exp\biggl(
      \frac{1-\lambda}{2\lambda} - \frac{R+C}{2\lambda^2n^2}
      + O(n^{-b})\biggr).\quad\qedsymbol
\]
\end{thm}

It is interesting to note that in the regular case $R=C=0$,
the average given in Theorem~\ref{permanent} is only higher than
Gurvits' lower bound~\cite{Gurvits} by a factor of 
$\lambda^{-1/2}(1+o(1))$.

\nicebreak
\section{Subdigraphs of random digraphs}\label{s:digraphs}

The adjacency matrix of a simple digraph is a square
$\{0,1\}$-matrix
with zero diagonal.  Therefore, Theorem~\ref{bigtheorem} can be
applied to enumerate digraphs with specified degrees, and the result
can then be used to explore the corresponding uniform probability space.

In this section, $\XX$ denotes a fixed simple digraph on the
vertices $\{w_1,\ldots,w_n\}$.
Let $\D(\svec,\tvec)$ be the
set of all simple digraphs on vertices $\{w_1,\ldots,w_n\}$
with out-degrees $\svec$ and in-degrees $\tvec$, and let
$\D(\svec,\tvec,\XX)$ be the subset of $\D(\svec,\tvec)$ containing
those digraphs which are arc-disjoint from~$\XX$.

For $1\leq j\leq n$ let
$x_j$, $y_j$ denote the out-degree and in-degree of
vertex $w_j$ in $\XX$, respectively.   The quantities
$s=t$, $\lambda$, $\delta_j$, $\eta_j$, $X$, $Y$, and so forth
are all defined by the same formulae as in Section~\ref{s:graphs}
with $m=n$.
In the definition of $Y$, the summation over $jk\in \XX$ should now
be interpreted as 
summation over $j,k\in \{ 1,\ldots, n\}$ such that $w_jw_k$
is an arc of $\XX$.
But note that $\lambda$ does not
represent the arc-density of a digraph in $\D(\svec,\tvec)$.
Instead the arc-density of a digraph in $\D(\svec,\tvec)$ is given
by 
\[ p = s/(n-1).\]

We begin with the basic enumeration result for digraphs.

\begin{thm}\label{bigdigraphtheorem}
For some $\eps > 0$, suppose that $s_j-s$, $x_j$, $t_j-s$
and $y_j$ are uniformly $O(n^{1/2+\eps})$ for $1\leq j, k\leq n$,
and $X= O(n^{1+2\eps})$, for $n\to\infty$.  Let $a,b>0$ be constants
such that $a+b< \nfrac{1}{2}$.  Suppose that $n\to\infty$ with
\[ \frac{(1-2\lambda)^2}{3A} \leq a\log n.\]
Then, provided $\eps > 0$ is small enough, we have
\begin{align*}
 & \card{\D(\svec,\tvec,\XX)}\\
 &{\quad}= \binom{n^2{-}X{-}n}{\lambda n^2}^{\!\!\!-1}
 \prod_{j=1}^n \binom{n{-}x_j{-}1}{s_j} \,\binom{n{-}y_j{-}1}{t_j} \\
  &{\qquad} \times 
   \exp\biggl( -\frac12\Bigl(1-\frac{R}{2An^2}\Bigr)
    \Bigl(1-\frac{C}{2An^2}\Bigr) 
    - \frac{Y +\sum_{j=1}^n (s_j-s)(t_j-s)}{2An^2} 
       + O(n^{-b})\biggr).
\end{align*}

\end{thm}

\begin{proof}
Let $\widetilde{\XX}$ be the bipartite graph
obtained from $\XX$ by replacing each vertex $w_j$ by two vertices
$u_j$, $v_j$, replacing each arc $w_jw_k$ of $\XX$
by the edge $u_jv_k$ of $\widetilde{\XX}$, and finally adding the 
perfect
matching $\{ u_jv_j\mid j=1,\ldots, n\}$ to the edge set of 
$\widetilde{\XX}$.
Then the degree sequences on the left and right of $\widetilde{\XX}$
and the total number of edges in $\widetilde{\XX}$ are given by
\[ \widetilde{x}_j = x_j + 1,\quad \widetilde{y}_j = y_j + 1,\quad
         \widetilde{X} = X + n,\]
         respectively ($1\leq j\leq n$).
The quantity $\widetilde{Y}$ for $\widetilde{\XX}$
satisfies
\[  \widetilde{Y} = Y + \sum_{j=1}^n (s_j-s)(t_j-s) + O(n^{2-b})\]
for any positive constant $b<1/2$.
Using this fact while applying Theorem~\ref{bigtheorem}
to $\widetilde{\XX}$ completes the proof.
\end{proof}

This formula for $\card{\D(\svec,\tvec,\XX)}$
has an intuitive interpretation
which is analogous to that given after Theorem~\ref{bigtheorem} for the
bipartite graph case.

Using this enumeration theorem, we can explore the uniform probability
space over $\D(\svec,\tvec)$. 
In each case, the proof is analogous to that of the corresponding
theorem for bipartite graphs in Section~\ref{s:graphs}.

\begin{thm}\label{avoidingXdig}
Under the conditions of Theorem~\ref{bigdigraphtheorem}, the following
are true for a random digraph $\BB\in\D(\svec,\tvec)$
if $\eps>0$ is small enough:
\begin{enumerate}
\itemsep=0pt
\item 
the probability that $\BB$ is arc-disjoint from $\XX$ is\/
$(1-p)^X \miss(n,n)$;
\item the probability $\BB$ contains $\XX$ as a subdigraph is\/
$p^X \hit(n,n)$,
\end{enumerate}
where $\miss(m,n)$ and $\hit(m,n)$ are defined in 
Theorem~\ref{avoidingXbip}.\quad\qedsymbol
\end{thm}

The special cases of $\miss(m,n)$ and $\hit(m,n)$ provided
by Corollaries~\ref{avoidingXflat} and~\ref{avoidingXreg} apply here
as well, as do the sufficient conditions~\eqref{biteme1}
and~\eqref{biteme2} for the
probabilities in Theorem~\ref{avoidingXdig} to asymptotically match
those in the standard random digraph model with arc probability~$p$.

Next suppose that each arc of $\XX$ has both ends in
$\{w_1,\ldots,w_J\}$.  Let $\XX_J$ be the subdigraph
of $\XX$ induced by those vertices.
For $\ell=1,2,3$, define
\begin{align*}
 \chi_\ell = \sum_{j=1}^J (x_j-p(J{-}1))^\ell,
\qquad
 \chi'_\ell = \sum_{k=1}^J (y_k-p(J{-}1))^\ell.
\end{align*}
Note that $\chi_1=\chi'_1=X-pJ(J-1)$.

\begin{thm}\label{inducedXdig}
Adopt the assumptions of Theorem~\ref{avoidingXdig}
with $s_j=s$ and $t_k=t$ for all~$j,k$, and
assume that $J=O(n^{1/2+\eps})$. 
The probability that a random digraph in $\D(\svec,\tvec)$
has $H_J$ as an induced subdigraph is
\begin{align*}
    p^X(1-p)^{J(J-1)-X} \exp\biggl(&
    \frac{J^2}{n} + \frac{(1-2\lambda)\chi_1}{2An}
           - \frac{\chi_1^2}{4An^2}\\
    &{} -\frac{(n+J)(\chi_2+\chi'_2)}{4An^2}
     -\frac{(1-2\lambda)(\chi_3+\chi'_3)}{24A^2n^2}
     + O(n^{-b})\biggr).\quad\qedsymbol
\end{align*}
\end{thm}

The argument of the exponential in Theorem~\ref{inducedXdig} is
$o(1)$ if $J^3=o(An)$.
So in that case, the probabilities of induced subdigraphs
asymptotically matches the standard random digraph model for arc
probability~$p$.

Let $\I(\XX)$ be the isomorphism class of $\XX$ and 
note that $|\I(\XX)| = n!/{\aut(\XX)}$, where
$\aut(\XX)$ is the number of automorphisms of~$\XX$.
By the same averaging technique as used to prove
Theorem~\ref{unlabelled}, we obtain the following.

\begin{thm}\label{unlabelleddig}
Suppose that the conditions of Theorem~\ref{bigdigraphtheorem}
apply with $x_j=y_j=x$ for all $j$.
Then the following is true of a random digraph $G$ in
$\D(\svec,\tvec)$:
\begin{enumerate}
\itemsep=0pt
\item the expected number of digraphs in $\I(\XX)$ that are
subgraphs of $G$ is
\[
   p^X\card{\I(\XX)}\exp\biggl(
      - \frac{(1-\lambda)\,x(x-2)}{2\lambda}
       - \frac{(R+C)x}{2\lambda^2n^2}
      + O(n^{-b})\biggr);
\]
\item the expected number of digraphs in $\I(\XX)$ that are
arc-disjoint from $G$ is
\[
   (1-p)^X\card{\I(\XX)}\exp\biggl(
      - \frac{\lambda x(x-2)}{2(1-\lambda)} 
       - \frac{(R+C)x}{2(1-\lambda)^2n^2}
      + O(n^{-b})\biggr). \quad\qedsymbol
\]
\end{enumerate}
\end{thm}

\nicebreak
\section{Proof of Theorem~\ref{bigtheorem}}\label{s:proof}

In the remainder of the paper we give the proof of
Theorem~\ref{bigtheorem}.
The overall method and many of the calculations will
parallel~\cite{CGM}, albeit with extra twists at each step,
so we acknowledge our considerable debt to Rod Canfield.

\begin{proof}[Outline of proof of Theorem~\ref{bigtheorem}.]
The basic idea is to identify $\card{\B(\svec,\tvec,\XX)}$ as a
coefficient in a multivariable generating function and to extract that
coefficient using the saddle-point method.
In Subsection~\ref{s:calculations},
we write $\card{\B(\svec,\tvec,\XX)}=P(\svec,\tvec,\XX)I(\svec,\tvec,\XX)$,
where $P(\svec,\tvec,\XX)$ is a rational expression and
$I(\svec,\tvec,\XX)$
is an integral in $m+n$ complex dimensions.  Both depend on the
location of the saddle point, which is the solution of some nonlinear
equations.  Those equations are solved in Subsection~\ref{s:radii},
and this leads to the value of $P(\svec,\tvec,\XX)$ in~\eqref{rad17}.
In~Subsections~\ref{s:integral}--\ref{s:complete},
the integral $I(\svec,\tvec,\XX)$ is estimated
in a small region $\R'$ defined in~\eqref{RSprime}.
The result is given by Lemma~\ref{Jintegral}
together with~\eqref{IvsJ}.  Finally, in Subsection~\ref{s:boxing},
it is shown that the integral $I(\svec,\tvec,\XX)$ restricted
to the exterior
of $\R'$ is negligible.  Theorem~\ref{bigtheorem} then follows
from~\eqref{start}, \eqref{rad17},
Lemmas~\ref{Jintegral}--\ref{boxing} and \eqref{IvsJ}.
\end{proof}

\medskip

We will use a shorthand notation for summation
over doubly subscripted variables.  If~$z_{jk}$ is a variable for
$1\leq j\leq m$ and $1\leq k\leq n$, then
\begin{align*}
  && z_{j\b} &= \sum_{k=1}^{n} z_{jk},&
   z_{\b k} &= \sum_{j=1}^{m} z_{jk}, &
   z_{\b\b} &= \sum_{j=1}^{m}\sum_{k=1}^{n} z_{jk}, &&\\
 && z_{j\c} &= \ksum z_{jk},&
   z_{\c k} &= \jsum z_{jk},&
   z_{\c\c} &= \jsum \ksum z_{jk}, &&
\end{align*}
for $1\leq j\leq m$ and $1\leq k\leq n$.

For $1\leq j\leq m$ and $1\leq k\leq n$ define $\h_{jk}=1$ 
if $u_jv_k$ is an edge of $\XX$ 
and $\h_{jk}=0$ otherwise. 
Then define the sets
\begin{align*}
 && X_j &= \{\, k \mid 1\le k\le n,\, \h_{jk}=1\,\}, &
 \barX_j &= \{\, k \mid 1\le k\le n,\, \h_{jk}=0\,\}, &&\\
 && Y_k &= \{\, j \mid 1\le j\le n,\, \h_{jk}=1\,\}, &
 \barY_k &= \{\, j \mid 1\le j\le n,\, \h_{jk}=0\,\}, &&
\end{align*}
for $1\le j\le m$ and $1\le k\le n$.

The notations $\sum_{jk\in\XX}$ and $\sum_{jk\in \barXX}$ indicate sums
over the sets $\{(j,k)\mid 1\le j\le m,\allowbreak 1\le k\le n, h_{jk}=1\}$
and $\{(j,k)\mid 1\le j\le m, 1\le k\le n, h_{jk}=0\}$, respectively,
and similarly for products.
We also define summations whose domain is limited by~$\XX$.
\begin{align*}
  && z_{j\b|\XX} &= \sum_{k\in X_j} z_{jk},&
   z_{\b k|\XX} &= \sum_{j\in Y_k} z_{jk}, &
   z_{\b\b|\XX} &= \sum_{jk\in\XX} z_{jk}, &&\\
  && z_{j\b|\barXX} &= \sum_{k\in \barX_j} z_{jk},&
   z_{\b k|\barXX} &= \sum_{j\in \barY_k} z_{jk}, &
   z_{\b\b|\barXX} &= \sum_{jk\in\barXX} z_{jk}. &&
\end{align*}

Under the assumptions of Theorem~\ref{bigtheorem}, we have
$m=\OO(n)$ and $n=\OO(m)$.  We also have that
$8\le A^{-1}\le O(\log n)$, so $A^{-1}=\OO(1)$.
More generally, $A^{c_1}m^{c_2+c_3\eps}n^{c_4+c_5\eps} =\OO(n^{c_2+c_4})$ 
if $c_1,c_2,c_3,c_4, c_5$ are constants.

We now show that the assumptions of Theorem~\ref{bigtheorem} imply that
\begin{equation}\label{oldassumptions}
m = o(A^2 n^{1+\eps}),\qquad n = o(A^2 m^{1+\eps}).
\end{equation}
If $A\geq \nfrac{3}{32}$ then \eqref{oldassumptions} follows immediately.
If $A < \nfrac{3}{32}$
then $(1-2\lambda)^2 > \nfrac{1}{4}$ and so the assumptions of
Theorem~\ref{bigtheorem} imply that $1/A = O\(\log n/(m/n+n/m)\)$.
This implies~\eqref{oldassumptions}.

\nicebreak
\subsection{Expressing the desired quantity as an integral}\label{s:calculations}

In this section we express $\card{\B(\svec,\tvec,\XX)}$
as a contour integral in
$(m+n)$-dimensional complex space, then begin to estimate its value
using the saddle-point method.

Firstly, notice that $\card{\B(\svec,\tvec,\XX)}$ is the coefficient of
$u_1^{s_1}\cdots u_m^{s_m}\,w_1^{t_1}\cdots w_n^{t_n}$ in the function
\[ \prod_{jk\in\barXX} \,(1+u_jw_k).\]
By Cauchy's coefficient theorem this equals
\[ \card{\B(\svec,\tvec,\XX)} = \frac{1}{(2\pi i)^{m+n}} \oint \cdots \oint
  \frac{\prod_{jk\in\barXX}
              (1+u_jw_k)}{u_1^{s_1+1}\cdots u_m^{s_m+1}
                 w_1^{t_1+1} \cdots w_n^{t_n+1}} \,
              du_1 \cdots du_m\, dw_1 \cdots dw_n,
\]
where each integral is along a simple closed contour enclosing the origin
anticlockwise.
It~will suffice to take each contour to be a circle;  specifically,
we will write
\[ u_j = q_j e^{i\theta_j} \quad \mbox{ and } \quad w_k = r_k e^{i\phi_k}\]
for $1\leq j\leq m$ and $1\leq k\leq n$.  Also define
\[ \lambda_{jk} = \frac{q_jr_k}{1+q_jr_k} \]
for $1\leq j\leq m$ and $1\leq k\leq n$.   Then
$\card{\B(\svec,\tvec,\XX)} = P(\svec,\tvec,\XX) I(\svec,\tvec,\XX)$
where
\begin{align}
\begin{split}\label{start}
 P(\svec,\tvec,\XX) &= \frac{\prod_{jk\in\barXX} (1 + q_jr_k)}
{(2\pi)^{m+n}\prod_{j=1}^m q_j^{s_j} \prod_{k=1}^n r_k^{t_k}}, \\
  I(\svec,\tvec,\XX) &=
  \int_{-\pi}^{\pi}\!\cdots \int_{-\pi}^\pi \frac{\prod_{jk\in\barXX}
    \(1 + \lambda_{jk}(e^{i(\theta_j+\phi_k)}-1)\)}
  {\exp(i\sum_{j=1}^m s_j\theta_j + i\sum_{k=1}^n t_k\phi_k)}
   \, d\thetavec d\phivec,
\end{split}
\end{align}
$\thetavec=(\t_1,\ldots,\t_m)$ and $\phivec=(\p_1,\ldots,\p_n)$.

We will choose the radii
$q_j$, $r_k$ so that there is no linear term in the logarithm
of the integrand of $I(\svec,\tvec,\XX)$ when expanded for
small $\thetavec,\phivec$.  This gives the equation
\[
 \sum_{jk\in\barXX}\lambda_{jk} (\theta_j+\phi_k)
 - \sum_{j=1}^m s_j\theta_j - \sum_{k=1}^n t_k\phi_k=0.
\]
For this to hold for all $\thetavec, \phivec$, we require
\begin{equation}\label{rad1}
\begin{split}
  \lambda_{j\b|\barXX} &= s_j \quad (1\le j\le m),\\
  \lambda_{\b k|\barXX} &= t_k \quad (1\le k\le n).
\end{split}
\end{equation}

The quantities $\lambda_{jk}$ have an interesting interpretation.
If edge $u_jv_k$ is chosen with probability $\lambda_{jk}$
independently for all $j,k\in\barXX$, then the expected degrees
are $\svec,\tvec$. 

In addition to the quantities defined before the statement of
Theorem~\ref{avoidingXbip}
we define for $j=1,\ldots, m$, $k=1,\ldots, n$,
\[
 \pTj_j = \sum_{k\in X_j}\dty_k, \quad
   \pSk_k = \sum_{j\in Y_k}\dsx_j.
\]

\nicebreak
\subsection{Locating the saddle-point}\label{s:radii}

In this subsection we solve \eqref{rad1} and derive some of the
consequences of the solution.  As with the whole paper, we
work under the assumptions of Theorem~\ref{bigtheorem}.

Change variables to $\{a_j\}_{j=1}^m$, $\{b_k\}_{k=1}^n$
as follows:
\begin{equation}\label{qrdef}
  q_j = r\frac{1+a_j}{1-r^2a_j},\quad
  r_k = r\frac{1+b_k}{1-r^2b_k},
\end{equation}
where
\[
  r= \sqrt{\frac{\lambda}{1-\lambda}}\;.
\]

Equation \eqref{rad1} is slightly underdetermined, which we will
exploit to impose an additional condition. If $\{q_j\}, \{r_k\}$
satisfy~\eqref{rad1} and $c>0$ is a constant, then $\{cq_j\},
\{r_k/c\}$ also satisfy~\eqref{rad1}.
{}From this we can see that, if there is a solution to~\eqref{rad1}
at all, there is one
for which $\sum_{j=1}^m (n-x_j)a_j<0$ and $\sum_{k=1}^n (m-y_k)b_k>0$, and
also a solution for which $\sum_{j=1}^m (n-x_j)a_j>0$ and
$\sum_{k=1}^n (m-y_k)b_k<0$.  It follows from the
Intermediate Value Theorem that there is a solution for which
\begin{equation}\label{rad4}
  \sum_{j=1}^m (n-x_j)a_j = \sum_{k=1}^n (m-y_k)b_k,
\end{equation}
so we will seek a common solution to~\eqref{rad1}
and \eqref{rad4}.

{}From \eqref{qrdef} we find that
\begin{equation}\label{rad2}
  \lambda_{jk}/\lambda = 1 + a_j + b_k + Z_{jk},
\end{equation}
where
\begin{equation}\label{Zjk}
 Z_{jk} = \frac{a_jb_k(1-r^2-r^2a_j - r^2b_k)}
               {1+r^2a_jb_k},
\end{equation}
and that equations~\eqref{rad1} can be rewritten as
\begin{equation}\label{rad3}
\begin{split}
  \frac{\dsx_j}{\lambda} &= (n-x_j)a_j
      + \sum_{k\in\barX_j}b_k + Z_{j\b|\barXX} \\
  \frac{\dty_k}{\lambda} &= (m-y_k)b_k
      + \sum_{j\in\barY_k}a_j + Z_{\b k|\barXX}\,.
\end{split}
\end{equation}

\noindent Summing \eqref{rad3} over all $j,k$, respectively,
we find in both cases that that
\begin{equation}\label{rad5}
 X= \sum_{j=1}^m (n-x_j)a_j + \sum_{k=1}^n (m-y_k)b_k
      + Z_{\b\b|\barXX}\,.
\end{equation}
Equations~\eqref{rad4} and~\eqref{rad5} together imply that
\[
   \sum_{j=1}^m (n-x_j)a_j = \sum_{k=1}^n (m-y_k)b_k 
    = \dfrac12 (X-Z_{\b\b|\barXX})\,.
\]
Substituting back into~\eqref{rad3}, we obtain
\begin{equation}\label{rad6}
\begin{split}
   a_j &= \mathbb{A}_j(a_1,\ldots,a_m,b_1,\ldots,b_n),\\
   b_k &= \mathbb{B}_k(a_1,\ldots,a_m,b_1,\ldots,b_n),
\end{split}
\end{equation}
for $1\le j\le m$, $1\le k\le n$, where
\begin{equation*}
\begin{split}
  \mathbb{A}_j(a_1,\ldots,a_m,b_1,\ldots,b_n) 
    &= \frac{\dsx_j}{\lambda n} - \frac{X}{2mn}
       + \frac{a_jx_j}{n} \\
     &{\kern1cm} - \frac{\sum_{k=1}^n y_kb_k}{mn}
       + \frac{\sum_{k\in X_j} b_k}{n}  
       - \frac{Z_{j\b|\barXX}}{n} 
       + \frac{Z_{\b\b|\barXX}}{2mn}, \\
  \mathbb{B}_k(a_1,\ldots,a_m,b_1,\ldots,b_n) 
    &= \frac{\dty_k}{\lambda m} - \frac{X}{2mn}
       + \frac{b_ky_k}{m} \\
     &{\kern1cm} - \frac{\sum_{j=1}^m x_ja_j}{mn}
       + \frac{\sum_{j\in Y_k} a_j}{m}  
       - \frac{Z_{\b k|\barXX}}{m} 
       + \frac{Z_{\b\b|\barXX}}{2mn}. \\
\end{split}
\end{equation*}

By the same argument as in~\cite{CGM},
equation \eqref{rad6} defines a convergent iteration
starting with $a_j=b_k=0$ for all~$j,k$.  Four iterations
give the following estimate of~$a_j$.  The value of $b_k$
follows by symmetry, while $Z_{jk}$ follows from~\eqref{Zjk}.
\begin{align*}
 a_j &= \frac{\dsx_j}{\lambda n}
         + \frac{\dsx_jx_j}{\lambda n^2}
         + \frac{\dsx_jx_j^2}{\lambda n^3}
         + \frac{\dsx_jx_j^3}{\lambda n^4}
         - \frac{X}{2mn}
         - \frac{(1-2\lambda)\dsx_j X}{4Amn^2}
         + \frac{\dsx_j^2 X}{4Amn^3} \\
     &{\quad} - \frac{x_jX}{mn^2}
         - \frac{x_j^2X}{mn^3}
         + \frac{\lambda(7-10\lambda)X^2}{16Am^2n^2}
         - \frac{3(1-2\lambda)\dsx_jx_jX}{4Amn^3}
         - \frac{(1-2\lambda)\pST}{4\lambda Am^2n^2} 
         \displaybreak[0]\\
     &{\quad} + \frac{\dsx_j C_2}{2\lambda Am^2n^2}
         + \frac{(1-2\lambda)\dsx_j^2 C_2}{4\lambda A^2m^2n^3}
         + \frac{\dsx_jx_jC_2}{\lambda Am^2n^3}
         - \frac{3XC_2}{8Am^3n^2}
         - \frac{XR_2}{8Am^2n^3} 
         \displaybreak[0]\\
     &{\quad} - \frac{(1-2\lambda)R_2C_2}{8\lambda A^2m^3n^3}
         - \frac{x_j\XS}{\lambda mn^3}
         - \frac{\YT}{\lambda m^2n}
         - \frac{(1-2\lambda)\dsx_j\YT}{2\lambda Am^2n^2}
         - \frac{x_j\YT}{\lambda m^2n^2}
         \displaybreak[0]\\
     &{\quad} - \frac{\XtwoS}{\lambda m^2n^2}
         - \frac{\VTone}{\lambda m^3n}
         + \frac{\pTj_j}{\lambda mn}
         + \frac{(1-2\lambda)\dsx_j\pTj_j}{2\lambda Amn^2}
         - \frac{\dsx_j^2\pTj_j}{2\lambda Amn^3}
         + \frac{x_j\pTj_j}{\lambda mn^2}
         \displaybreak[0]\\
     &{\quad} + \frac{x_j^2\pTj_j}{\lambda mn^3}
         + \frac{(1-2\lambda)\dsx_jx_j\pTj_j}{\lambda Amn^3}
         - \frac{3(1-2\lambda)X\pTj_j}{4Am^2n^2}
         + \frac{(1-2\lambda)\pTj_j^2}{2\lambda Am^2n^2}
         \displaybreak[0]\\
     &{\quad} + \Bigl(\frac{R_2}{n}+\frac{C_2}{m}\Bigr)
               \frac{\pTj_j}{2\lambda Am^2n^2}
         + \frac{1}{\lambda m^2n^2}\sum_{(j',k')}\dsx_{j'}y_{k'}
         + \frac{1}{\lambda m^3n}\sum_{k\in X_j}\dty_{k}y_{k}^2
         \displaybreak[0]\\
     &{\quad} - \frac{X}{m^2n^2}\sum_{k\in X_j}y_k
         - \frac{\dsx_j}{2\lambda Am^2n^2}\sum_{k\in X_j}\dty_k^2
         + \frac{(1-2\lambda)}{2\lambda Am^2n^2}
                \sum_{(j',k')}\dsx_{j'}\dty_{k'}
         \displaybreak[0]\\
     &{\quad} + \Bigl(\frac{1}{\lambda m^2n^2} 
                 + \frac{(1-2\lambda)\dsx_j}{2\lambda Am^2n^3}
                 + \frac{x_j}{\lambda m^2n^3}\Bigr)
                   \Bigl( n\sum_{k\in X_j} \dty_ky_k
                         + m\sum_{(j',k')}\dsx_{j'}\Bigr)
         \displaybreak[0]\\
     &{\quad} + \frac{1}{\lambda mn^3}\sum_{(j',k')}\dsx_{j'} x_{j'}
         + \frac{1}{\lambda m^2n^2}
                      \sum_{(j',k')}\sum_{k''\in X_{j'}}\dty_{k''}
         + \OO(n^{-5/2}),
\end{align*}
where the notation $\sum_{(j',k')}$ means
$\sum_{k'\in X_j}\sum_{j'\in Y_{k'}}$.

A sufficient approximation of $\lambda_{jk}$ is given by substituting
this estimate into~\eqref{rad2}.   In~evaluating
the integral $I(\svec,\tvec,\XX)$, the following approximations
will be required:
\begin{gather}
\begin{split}\label{rad10}
  \lambda_{jk}(1-\lambda_{jk}) &= \lambda(1-\lambda)
        + \frac{(1-2\lambda)\dsx_j}{n} + \frac{(1-2\lambda)\dty_k}{m}
         - \frac{\dsx_j^2}{n^2} - \frac{\dty_k^2}{m^2} \\
        &\quad{}
        + \frac{(1-12A)\dsx_j\dty_k}{2Amn}
        + \frac{(1-2\lambda)\dsx_jx_j}{n^2}
        + \frac{(1-2\lambda)\dty_ky_k}{m^2} \\
        &\quad{} + \frac{(1-2\lambda)(\pTj_j+\pSk_k-\lambda X)}{mn}
                + \OO(n^{-3/2}),
\end{split} \displaybreak[0]\\[0.5ex]
\begin{split}\label{rad11}
  \lambda_{jk}(1-\lambda_{jk})(1-2\lambda_{jk})
 &= \lambda(1-\lambda)(1-2\lambda)
        + \frac{(1-12A)\dsx_j}{n}\\
        &\qquad{} + \frac{(1-12A)\dty_k}{m} + \OO(n^{-1}),
\end{split} \\[0.5ex]
\begin{split}\label{rad12}
  \lambda_{jk}(1-\lambda_{jk})(1-6\lambda_{jk}+6\lambda_{jk}^2)
     &= \lambda(1-\lambda)(1-12A) + \OO(n^{-1/2}).
\end{split}
\end{gather}

We now estimate the factor $P(\svec,\tvec,\XX)$.
If
\[\Lambda = \prod_{jk\in\barXX}
          \lambda_{jk}^{\lambda_{jk}} (1-\lambda_{jk})^{1-\lambda_{jk}}\]
then
\begin{align*}
\Lambda^{-1} &= \prod_{jk\in\barXX}
  \biggl(\frac{1 + q_jr_k}{q_j r_k}\biggr)^{\!\!\lambda_{jk}}
        (1 + q_j r_k)^{1-\lambda_{jk}}\\
 &= \prod_{jk\in\barXX} \,(1 + q_j r_k) \,
  \biggl(\,\prod_{j=1}^m q_j^{\lambda_{j\b|\barXX}}\,
       \prod_{k=1}^n r_k^{ \lambda_{\b k|\barXX}}\biggr)^{\!\!-1}\\
 &= \prod_{jk\in\barXX} \,(1+q_j r_k)\,
  \prod_{j=1}^m q_j^{-s_j} \prod_{k=1}^n r_k^{-t_k}
\end{align*}
using \eqref{rad1}.  Therefore, the factor
$P(\svec,\tvec,\XX)$ in front
of the integral in \eqref{start} is given by
\[ P(\svec,\tvec,\XX) = (2\pi)^{-(m+n)} \, \Lambda^{-1}.\]
We proceed to estimate $\Lambda$.
Writing $\lambda_{jk}=\lambda(1+z_{jk})$, we have
\begin{equation}\label{rad13}
\begin{split}
  &\kern-2mm\log\biggl( \frac{\lambda_{jk}^{\lambda_{jk}}
     (1-\lambda_{jk})^{1-\lambda_{jk}}}
                 {\lambda^\lambda(1-\lambda)^{1-\lambda}}\biggr)
     = \lambda z_{jk}\log\biggl(\frac{\lambda}{1-\lambda}\biggr)\\
      &{\qquad} +
     \frac{\lambda}{2(1-\lambda)}z_{jk}^2 - \frac{\lambda(1-2\lambda)}
          {6(1-\lambda)^2} z_{jk}^3
       + \frac{\lambda(1-3\lambda+3\lambda^2)}{12(1-\lambda)^3} z_{jk}^4
       + O\biggl(\frac{z_{jk}^5}{(1-\lambda)^4}\biggr).
\end{split}
\end{equation}
We know from \eqref{rad1} that $\lambda_{\b\b|\barXX} = \lambda mn$,
which implies that
$z_{\b\b|\barXX}=X$, hence the first term on the right side
of~\eqref{rad13}
contributes $\lambda^{\lambda X}(1-\lambda)^{-\lambda X}$ to $\Lambda$.
Now using \eqref{rad2}
we can write $z_{jk} = a_j + b_k + Z_{jk}$ and apply the above
estimates to obtain
\begin{align}
  \Lambda &= \( \lambda^\lambda(1-\lambda)^{1-\lambda} \)^{mn}
   (1-\lambda)^{-X} \notag\\
  &{\kern8mm}\times 
   \exp\biggr( \frac{R_2}{4An} + \frac{C_2}{4Am} 
        + \frac{R_2C_2}{8A^2m^2n^2}
        - \frac{\lambda^2 X^2}{4Amn}
        - \frac{(1-2\lambda)}{24A^2}
          \Bigl(\frac{R_3}{n^2}+\frac{C_3}{m^2}\Bigr) \label{rad14}\\
  &{\kern22mm}  + \frac{(1-6A)}{96A^3}
               \Bigl(\frac{R_4}{n^3}+\frac{C_4}{m^3}\Bigr)
           + \frac{\pST}{2Amn}
    + \frac{\XStwo}{4An^2} + \frac{\YTtwo}{4Am^2}
    + \OO(n^{-1/2}) \biggl). \notag
\end{align}

As in \cite{CGM}, our answer will be simpler when written
in terms of binomial coefficients. 
Using an accurate approximation of the binomial coefficients
(such as~\cite[Equation 18]{CGM}), we obtain that
\begin{equation}\label{rad16}
\begin{split}
\binom{mn{-}X}{\lambda mn}^{\!\!-1}
  &\, \prod_{j=1}^m\binom{n{-}x_j}{s_j}\prod_{k=1}^n\binom{m{-}y_k}{t_k}
  = \frac{(\lambda^\lambda(1-\lambda)^{1-\lambda})^{-mn}(1-\lambda)^X}
          {(4\pi A)^{(m+n-1)/2} m^{(n-1)/2} n^{(m-1)/2}}\\
  &{}\times \exp\biggl( - \frac{R_2}{4An} - \frac{C_2}{4Am}
    - \frac{1-2A}{24A}\Bigl(\frac{m}{n}+\frac{n}{m}\Bigr)
   + \frac{1-4A}{16A^2}\Bigl(\frac{R_2}{n^2}+\frac{C_2}{m^2}\Bigr)\\
  &{\kern14mm}
  + \frac{1-2\lambda}{24A^2}\Bigl(\frac{R_3}{n^2}+\frac{C_3}{m^2}\Bigr)
   - \frac{1-6A}{96A^3}\Bigl(\frac{R_4}{n^3}+\frac{C_4}{m^3}\Bigr) \\
  &{\kern14mm}
   + \frac{\lambda^2 X(m+n+X)}{4Amn}
   - \frac{\XStwo}{4An^2}
   - \frac{\YTtwo}{4Am^2}
   + \OO(n^{-1/2})\biggl).
\end{split}
\end{equation}
Putting \eqref{rad14} and \eqref{rad16} together, we find that
\begin{equation}\label{rad17}
\begin{split}
P(\svec,\tvec,\XX)
   &= \Lambda^{-1}(2\pi)^{-(m+n)}\\
   &= \frac{A^{(m+n-1)/2} m^{(n-1)/2} n^{(m-1)/2}}
                       {2\pi^{(m+n+1)/2}}
      \;\binom{mn{-}X}{\lambda mn}^{\!\!-1}
      \prod_{j=1}^{m}\binom{n{-}x_j}{s_j}\prod_{k=1}^{n}\binom{m{-}y_k}{t_k}  \\
   & {\qquad} \times
     \exp\biggl( \frac{1-2A}{24A}\Bigl(\frac{m}{n}+\frac{n}{m}\Bigr)
      - \frac{R_2C_2}{8A^2m^2n^2}
      - \frac{1-4A}{16A^2}\Bigl(\frac{R_2}{n^2}+\frac{C_2}{m^2}\Bigr) \\
   & {\kern23mm}
      - \frac{\lambda^2 X}{4A}\Bigl(\frac{1}{m}+\frac{1}{n}\Bigr)
      - \frac{\pST}{2Amn}
      + \OO(n^{-1/2})\biggr).
\end{split}
\end{equation}

\nicebreak
\subsection{Evaluating the integral}\label{s:integral}

Our next task is to evaluate the integral $I(\svec,\tvec,\XX)$
given~\eqref{start}.

Let $C$ be the ring of real numbers modulo $2\pi$, which we can interpret
as points on a circle, and let $z$ be the canonical
mapping from $C$ to the real interval $(-\pi, \pi]$.  An \emph{open
half-circle} is $C_t = (t-\pi/2, t+\pi/2) \subseteq C$ for some $t$.
Now define
\[ \Chat^N = \{\, \vvec = (v_1,\ldots, v_N) \in C^N \mid
  v_1, \ldots, v_N \in C_t \mbox{ for some } t\in\mathbb{R}\,\}.\]
If $\vvec = (v_1,\ldots, v_N)\in C_0^N$ then define
\[ \avv = z^{-1}\biggl( \frac{1}{N}
  \sum_{j=1}^N z(v_j)\biggr).\]
More generally, if $\vvec\in C_t^N$ then define
$\avv = t + \overline{(v_1-t,\ldots, v_N-t)}$.
The function $\vvec\to \avv$ is
well-defined and continuous for $\vvec\in\Chat^N$.

Let $\R$
denote the set of vector pairs $(\thetavec ,\phivec)
\in \Chat^m\times\Chat^n$ such that
\begin{align}\label{Rdef}
\begin{split}
\abs{\avtheta + \avphi}
  &\leq  (mn)^{-1/2 + 2\eps},\\
\abs{\that_j} &\leq n^{-1/2 + \eps} \qquad (1\leq j\leq m),\\
\abs{\phat_k} &\leq m^{-1/2 + \eps} \qquad (1\leq k\leq n),
\end{split}
\end{align}
where $\that_j = \theta_j - \avtheta$
and $\phat_k = \phi_k - \avphi$.  In
this definition, values are considered in $C$.  The constant
$\eps$ is the sufficiently-small value required by
Theorem~\ref{bigtheorem}.

\smallskip

Let $I_{\R''}(\svec,\tvec,\XX)$ denote the integral
$I(\svec,\tvec,\XX)$ restricted to any region $\R''$.
In this subsection, we estimate $I_{\R'}(\svec,\tvec,\XX)$
in a certain region $\R'\supseteq\R$.
In Subsection~\ref{s:boxing} we will show that
the remaining parts of $I(\svec,\tvec,\XX)$ are negligible.
We begin by analysing
the integrand in $\R$, but for future use when we
expand the region to~$\R'$ (to be defined in~\eqref{RSprime}),
note that all the approximations
we establish for the integrand in $\R$ also hold in the
superset of $\R'$ defined~by
\begin{align}\begin{split}\label{bigR}
\abs{\avtheta + \avphi}
  &\leq  3(mn)^{-1/2 + 2\eps},\\
\abs{\that_j} &\leq 3n^{-1/2 + \eps} \quad (1\leq j\leq m-1),\\
\abs{\that_m} &\leq 2n^{-1/2 + 3\eps},\\
\abs{\phat_k} &\leq 3m^{-1/2 + \eps} \quad (1\leq k\leq n-1),\\
\abs{\phat_n} &\leq 2m^{-1/2 + 3\eps}.
\end{split}\end{align}

Define $\thetahatvec=(\that_1,\ldots,\that_{m-1})$ and
$\phihatvec=(\phat_1,\ldots,\phat_{n-1})$.
Let $T_1$ be the transformation
$T_1(\thetahatvec,\phihatvec,\nu,\psi)=(\thetavec,\phivec)$
defined by
\[ \nu = \avtheta + \avphi,
  \quad  \psi = \avtheta -
\avphi,\]
together with $\that_j=\t_j-\avtheta$
($1\leq j\leq m-1$) and $\phat_k=\p_k-\avphi$ ($1\leq k\leq n-1$).
We also define the 1-many transformation $T_1^*$ by
\[
    T_1^*(\thetahatvec,\phihatvec,\nu) =
    \bigcup_\psi \, T_1(\thetahatvec,\phihatvec,\nu,\psi).
\]

After applying the transformation $T_1$ to $I_\R(\svec,\tvec,\XX)$,
the new integrand is easily seen to be independent of
$\psi$, so we can multiply by the range of $\psi$ and remove
it as an independent variable.  Therefore, we can continue with
an $(m+n-1)$-dimensional integral over $\S$ such that
$\R=T_1^*(\S)$.   More generally, if
$\S''\subseteq (-\tfrac12\pi,\tfrac12\pi)^{m+n-2}\times(-2\pi,2\pi]$
and $\R''=T_1^*(\S'')$, we have
\begin{equation}\label{IvsJ}
I_{\R''}(\svec,\tvec,\XX) =  2\pi m n \int_{\S''} 
   G(\thetahatvec,\phihatvec,\nu)
     \, d\thetahatvec  d\phihatvec d\nu,
\end{equation}
where $G(\thetahatvec,\phihatvec,\nu)
   =F\(T_1(\thetahatvec,\phihatvec,\nu,0)\)$
with $F(\thetavec,\phivec)$ being the integrand of $I(\svec,\tvec,\XX)$.
The factor $2\pi mn$ combines the range of $\psi$,
which is $4\pi$, and the Jacobian of $T_1$, which is~$mn/2$.

Note that  $\S$ is defined by
the same inequalities \eqref{Rdef} as define~$\R$.  The first inequality
is now $\abs{\nu}\leq (mn)^{-1/2 + 2\eps}$ and the bounds on
\[ \that_m = - \jsum \that_j~\text{ and }~
    \phat_m = - \ksum \phat_k
\]
still apply even though these are no
longer variables of integration.

\medskip

In the region $\S$, the integrand of \eqref{IvsJ} can be expanded as
\begin{align*}
 G(\thetahatvec,\phihatvec,\nu) &= \exp\biggl( -\sum_{jk\in\barXX}
        (A+\alpha_{jk})(\nu+\that_j+\phat_k)^2
           - i\sum_{jk\in\barXX}(A_3+\beta_{jk})(\nu+\that_j+\phat_k)^3 \\
     & \qquad\quad{}
             + \sum_{jk\in\barXX}(A_4+\gamma_{jk})(\nu+\that_j+\phat_k)^4
             + O\Bigl(A\,\sum_{jk\in\barXX}\,\abs{\nu+\that_j+\phat_k}^5\,\Bigr)
              \biggr) \\
  &= \exp\biggl( -\sumjkmn
        (A+\alpha_{jk})(\nu+\that_j+\phat_k)^2
                  - i\sumjkmn(A_3+\beta_{jk})(\nu+\that_j+\phat_k)^3 \\
     & \qquad\quad{}
                  + \sumjkmn(A_4+\gamma_{jk})(\nu+\that_j+\phat_k)^4
                  + \sum_{jk\in\XX} A(\nu+\that_j+\phat_k)^2
                  + \OO(n^{-1/2})
              \biggr).
\end{align*}
Here $\alpha_{jk}$, $\beta_{jk}$, and $\gamma_{jk}$ are defined by
\begin{align}\label{AlBetGamDef}
\begin{split}
\dfrac{1}{2}\lambda_{jk}(1-\lambda_{jk}) &=  A + \alpha_{jk}, \\
\dfrac{1}{6}\lambda_{jk}(1-\lambda_{jk})(1-2\lambda_{jk})
   &= A_3 + \beta_{jk}, \\
\dfrac{1}{24}\lambda_{jk}(1-\lambda_{jk})(1-6\lambda_{jk}+6\lambda_{jk}^2)
   &= A_4 + \gamma_{jk},
\end{split}
\end{align}
where
\[ A = \dfrac{1}{2}\lambda(1-\lambda),~
 A_3 = \dfrac{1}{6}\lambda(1-\lambda)(1-2\lambda), \text{~and~}
  A_4 = \dfrac{1}{24}\lambda(1-\lambda)(1-6\lambda + 6\lambda^2).
\]
Approximations for $\alpha_{jk}$, $\beta_{jk}$, $\gamma_{jk}$ were
given in \eqref{rad10}--\eqref{rad12}.  Note that $\alpha_{jk}$
in this paper is slightly different from in~\cite{CGM}, but it
is still true that
$\alpha_{jk}, \beta_{jk}, \gamma_{jk}=\OO(n^{-1/2})$ uniformly
over~$j,k$.

\nicebreak
\subsection{Another change of variables}\label{s:change}

We now make a second change of variables
$(\thetahatvec,\phihatvec,\nu) =T_2(\zetavec,\xivec,\nu)$,
where $\zetavec=(\zeta_1,\ldots,\zeta_{m-1})$ and
$\xivec=(\xi_1,\ldots,\xi_{n-1})$,
whose purpose is to almost diagonalize the quadratic part
of~$G$.  The diagonalization will be completed in the
next subsection.
The transformation $T_2$ is defined as follows.
For $1\leq j\leq m-1$ and $1\leq k\leq n-1$ let
\[ \that_j = \zeta_j + c\pi_1,\quad \phat_k = \xi_k + d\rho_1, \]
where
\[ c=-\frac{1}{m+m^{1/2}}
   \quad \mbox{and} \quad d= - \frac{1}{n+n^{1/2}}\]
and, for $1\leq h\leq 4$,
\[ \pi_h = \jsum \zeta_j^h,\quad \rho_h = \ksum
           \xi_k^h.\]
The Jacobian of the transformation is $(mn)^{-1/2}$.
In \cite{CM}, this transformation was seen to exactly diagonalize
the quadratic part of the integrand in the semiregular case.
In the present irregular case, the diagonalization is no longer
exact but still provides useful progress.

By summing the equations $\that_j = \zeta_j + c\pi_1$ and
$\phat_k = \xi_k + d\rho_1$, we find that
\begin{align}\label{pi1bound}
\begin{split}
   \pi_1 = m^{1/2}\jsum\that_j,\quad \abs{\pi_1}\le m^{1/2}n^{-1/2+\eps}, \\
   \rho_1 = n^{1/2}\ksum\phat_k,\quad \abs{\rho_1}\le n^{1/2}m^{-1/2+\eps},
\end{split}
\end{align}
where the inequalities come from the bounds on $\that_m$ and $\phat_n$.
This implies that
\begin{align*}
   \zeta_j &= \that_j+\OO(n^{-1})\quad(1\le j\le m-1), \\
   \xi_k &= \phat_k+\OO(n^{-1})\quad(1\le k\le n-1).
\end{align*}
The transformed region of integration is $T_2^{-1}(\S)$, but
for convenience we will expand it a little to be the region
defined by the inequalities
\begin{align}\label{zetaxibox}
\begin{split}
   \abs{\zeta_j} &\le \tfrac32 n^{-1/2+\eps}\quad(1\le j\le m-1), \\
   \abs{\xi_k} &\le \tfrac32 m^{-1/2+\eps}\quad(1\le k\le n-1), \\
   \abs{\pi_1} &\le m^{1/2}n^{-1/2+\eps}, \\
   \abs{\rho_1} &\le n^{1/2}m^{-1/2+\eps}, \\
   \abs{\nu} &\le (mn)^{-1/2+2\eps}\,.
\end{split}
\end{align}

We now consider the new integrand $E_1=\exp(L_1) = G\circ T_2$.
As in~\cite{CM}, the semiregular parts of $L_1$ (those not
involving $\alpha_{jk}$, $\beta_{jk}$, $\gamma_{jk}$ or~$\XX$)
transform to 
\begin{align}\label{regbits}
\begin{split}
&-Amn\nu^2 - A n \pi_2 - A m \rho_2
   -3i A_3 n\nu \pi_2 - 3iA_3 m\nu \rho_2
  + 6 A_4\pi_2\rho_2 \\
 &\quad  {} - i A_3 n\pi_3
        - i A_3 n\rho_3
- 3i A_3 cn \pi_1\pi_2
- 3iA_3 dm \rho_1\rho_2
 + A_4 n\pi_4 + A_4 m\rho_4 + \OO(n^{-1/2}).
\end{split}
\end{align}
To see the effect of the transformation on the irregular parts of the
integrand, write $\zeta_m=\that_m-c\pi_1$ and
$\xi_n=\that_n-d\rho_1$.  {}From~\eqref{pi1bound} we can see
that $\zeta_m=\OO(n^{-1/2})$ and $\xi_n=\OO(n^{-1/2})$.
Thus we have, for all $1\le j\le m$ and $1\le k\le n$,
$\zeta_j+\xi_k = \OO(n^{-1/2})$ and $c\pi_1+d\rho_1+\nu=\OO(n^{-1})$.
Recalling also that $\alpha_{jk}, \beta_{jk}, \gamma_{jk} = \OO(n^{-1/2})$,
we have
\begin{align*}
  \sum_{j=1}^m\sum_{k=1}^n\alpha_{jk} (\nu+\that_j+\phat_k)^2\\[-2ex]
  &\kern-2cm{}= \sum_{j=1}^m\sum_{k=1}^n\alpha_{jk} 
         \((\zeta_j+\xi_k)^2+2(\zeta_j+\xi_k) (\nu+c\pi_1+d\rho_1)\) + \OO(n^{-1/2}), \\
  \sum_{j=1}^m\sum_{k=1}^n\beta_{jk} (\nu+\that_j+\phat_k)^3
  &= \sum_{j=1}^m\sum_{k=1}^n\beta_{jk}(\zeta_j+\xi_k)^3 + \OO(n^{-1/2}), \\
  \sum_{j=1}^m\sum_{k=1}^n\gamma_{jk} (\nu+\that_j+\phat_k)^4
  &= \OO(n^{-1/2}), \\
  \sum_{jk\in\XX} (\nu+\that_j+\phat_k)^2
  &= \sum_{jk\in\XX} (\zeta_j+\xi_k)^2 + \OO(n^{-1/2})
\end{align*}
Moreover, the terms on the right sides of the above that involve
$\zeta_m$ or $\xi_n$ contribute only $\OO(n^{-1/2})$ in total, so we can drop them.
Combining this with~\eqref{regbits}, we have
\begin{equation}\label{Helper}
\begin{split}
L_1 &=  -Amn\nu^2 - A n \pi_2 - A m \rho_2
   -3i A_3 n\nu \pi_2 - 3iA_3 m\nu \rho_2
  + 6 A_4\pi_2\rho_2 \\[0.6ex]
 & \hspace*{5mm} {} - i A_3 n\pi_3
        - i A_3 n\rho_3
- 3i A_3 cn \pi_1\pi_2
- 3iA_3 dm \rho_1\rho_2
 + A_4 n\pi_4 + A_4 m\rho_4 \\
 & \hspace*{5mm} {}
   - \jsum\ksum
   \alpha_{jk}\((\zeta_j+\xi_k)^2+2(\zeta_j+\xi_k) (\nu+c\pi_1+d\rho_1)\)
          \\
 & \hspace*{5mm} {}
   - i \jsum\ksum \beta_{jk}(\zeta_j+\xi_k)^3
   + A\sum_{jk\in\XX} (\zeta_j+\xi_k)^2
   + \OO(n^{-1/2}).
\end{split}
\end{equation}

\nicebreak
\subsection{Completing the diagonalization}\label{s:diagonalize}

The quadratic form in $L_1$ is the following function of the $m+n-1$
variables $\zetavec,\xivec,\nu$:
\begin{align}\begin{split}\label{Qvalue}
Q &= -Amn\nu^2 - An \pi_2 - Am \rho_2 
      + A\sum_{jk\in\XX} (\zeta_j+\xi_k)^2 \\
    &\quad {} - \jsum\ksum
\alpha_{jk}\((\zeta_j+\xi_k)^2 + 2(\zeta_j+\xi_k)(\nu+c\pi_1+d\rho_1)\).
\end{split}\end{align}
We will make a third change of variables,
$(\zetavec,\xivec,\nu)=T_3 (\sigmavec,\tauvec,\mu)$, that
diagonalizes this quadratic form,
where $\sigmavec=(\sigma_1,\ldots,\sigma_{m-1})$
and $\tauvec=(\tau_1,\ldots,\tau_{n-1})$.
This is achieved using a slight extension
of \cite[Lemma 3.2]{MWtournament}.

\begin{lemma}
\label{diagonalize}
Let $\UU$ and $\YY$ be square matrices of the same order, such that
$\UU^{-1}$ exists and all the eigenvalues of
$\UU^{-1} \YY$ are less than 1 in absolute value.
Then
\[ (\II + \YY\UU^{-1})^{-1/2} \,(\UU+\YY)\,
  (\II+\UU^{-1}\YY)^{-1/2} = \UU, \]
where the fractional powers are defined by the binomial expansion.
\quad\qedsymbol
\end{lemma}

Note that $\UU^{-1}\YY$ and $\YY\UU^{-1}$ have the same eigenvalues, so
the eigenvalue condition on $\UU^{-1}\YY$ applies equally to $\YY\UU^{-1}$.
If we also have that both $\UU$ and $\YY$ are symmetric, then
$(\II+\YY\UU^{-1})^{-1/2}$ is the transpose of
$(\II + \UU^{-1}\YY)^{-1/2}$, as proved in~\cite{CGM}.
Let $\VV$ be the symmetric matrix associated with the quadratic form~$Q$.
 Write $\VV=\Vd + \Vnd$ where $\Vd$ has all off-diagonal
entries equal to zero and matches $V$ on the diagonal entries,
and $\Vnd$ has all diagonal entries zero and matches $\VV$ on the
off-diagonal entries.
We will apply Lemma~\ref{diagonalize}
with $\UU = \Vd$ and $\YY=\Vnd$.
Note that $\Vd$ is invertible and that both $\Vd$
and $\Vnd$ are symmetric.
Let $T_3$ be the transformation given by
$T_3(\sigmavec,\tauvec,\mu)^T = (\zetavec,\xivec,\nu)^T=
(\II+\Vdinv \Vnd)^{-1/2}(\sigmavec,\tauvec,\mu)^T$.
If the eigenvalue condition of Lemma~\ref{diagonalize} is satisfied then
this transformation diagonalizes the quadratic form $Q$, keeping
the diagonal entries unchanged.

{}From the formula for $Q$ we extract the following coefficients,
which tell us the diagonal and off-diagonal entries of $\VV$.
Define $x'_j = x_j-\h_{jn}$ for $1\le j\le m-1$, and
$y'_k=y_k-\h_{mk}$ for $1\le k\le n-1$.  Then:
\begin{align*}
[\zeta_j^2]\,Q &= -An-(1+2c)\alpha_{j\c}+Ax'_j,\\
[\xi_k^2]\,Q &= -Am-(1+2d)\alpha_{\c k}+Ay'_k,\\
[\nu^2]\,Q &= -Amn, 
\displaybreak[0]\\
[\zeta_{j_1}\zeta_{j_2}]\,Q &= -2c(\alpha_{j_1\c}+\alpha_{j_2\c})
 \kern20mm(j_1\ne j_2),\\
[\zeta_j \xi_k]\,Q 
 &= -2\alpha_{jk} - 2d\alpha_{j\c} - 2c\alpha_{\c k} + 2A\h_{jk},\\
[\xi_{k_1}\xi_{k_2}]\,Q &= -2d(\alpha_{\c k_1}+\alpha_{\c k_2})
 \kern20mm(k_1\ne k_2),\\
[\zeta_j \nu]\,Q &= -2\alpha_{j\c}, \\
[\xi_k \nu]\,Q &= -2\alpha_{\c k}.
\end{align*}
Using these equations we find that all off-diagonal entries of
$\Vdinv\Vnd$ are $\OO(n^{-3/2})$, except for the column
corresponding to $\nu$, which has off-diagonal entries of
size $\OO(n^{-1/2})$, and the entries corresponding to $\zeta_j\xi_k$
for $\h_{jk}=1$, which have size~$\OO(n^{-1})$.
Similarly, the off-diagonal entries of $\Vnd\Vdinv$
are all $\OO(n^{-3/2})$, except for the row corresponding to $\nu$,
which has off-diagonal entries of size $\OO(n^{-1/2})$, and
the entries corresponding to $\zeta_j\xi_k$
for $\h_{jk}=1$, which have size~$\OO(n^{-1})$.  To see that
these conditions imply that the eigenvalues of
$\Vdinv\Vnd$
are less than one, recall that the value of any matrix norm is
greater than or equal to the greatest absolute value of an
eigenvalue.  The $\infty$-norm (maximum row sum of
absolute values) of $\Vdinv\Vnd$ is $\OO(n^{-1/2})$, so the
eigenvalues are all~$\OO(n^{-1/2})$.

We also need to know the Jacobian of the transformation $T_3$.

\begin{lemma}[\cite{CGM}]\label{l:detexp}
  Let $\MM$ be a matrix of order $O(m+n)$ with all eigenvalues
  uniformly $\OO(n^{-1/2})$.  Then
  \[
    \det(\II+\MM) = \exp\(\tr\MM - \dfrac12\tr\MM^2 + \OO(n^{-1/2})\).
    \quad\qedsymbol
  \]
\end{lemma}

Let $\MM=\Vdinv\Vnd$.
As noted before, the eigenvalues of $\MM$ are all $\OO(n^{-1/2})$ so
Lemma~\ref{l:detexp} applies.
Noting that
$\tr(\MM)=0$ and calculating that $\tr(\MM^2) = \OO(n^{-1})$, we
conclude that the Jacobian of $T_3$ is
\[
  \det\((\II+\MM)^{-1/2}\) = \(\det(\II+\MM)\)^{-1/2}
 = 1+ \OO(n^{-1/2}). 
\]
% [[BDM: Actually it is $1+\OO(n^{-1})$ but we don't need that.]]

\medskip

To derive $T_3$ explicitly, we can expand
$(\II+\Vdinv\Vnd)^{-1/2}$ while noting that
$\alpha_{j\c} = O(n^{1/2+\eps})$ for all~$j$,
$\alpha_{\c k} = O(m^{1/2 + \eps})$ for all~$k$,
$\alpha_{\c\c} = O(mn^{2\eps} + nm^{2\eps})$,
$R\le mn^{1+2\eps}$ and $C\le nm^{1+2\eps}$.

\nicebreak
This gives
\begin{align*}
\sigma_j &= \zeta_j + \sum_{j'=1}^{m-1} \Bigl(\frac{c(\alpha_{j\c}
    + \alpha_{j'\c})}{2An} + \OO(n^{-2})\Bigr) \zeta_{j'} \\
    &{\qquad}+
        \ksum \Bigl(\frac{\alpha_{jk} + d\alpha_{j\c} +
          c\alpha_{\c k}}{2An} + \OO(n^{-2})\Bigr) \xi_k
   + \Bigl(\frac{\alpha_{j\c}}{2An} + \OO(n^{-1})\Bigr)\nu + \OO(n^{-3/2}),
       \displaybreak[0]\\
\tau_k &= \xi_k + \jsum \Bigl(\frac{\alpha_{jk} +
      d\alpha_{j\c} + c\alpha_{\c k}}{2Am} 
        + \OO(n^{-2})\Bigr)\zeta_j \\
    &{\qquad}
    + \sum_{k'=1}^{n-1} \Bigl(\frac{d(\alpha_{\c k} + \alpha_{\c k'})}
  {2Am} + \OO(n^{-2})\Bigr)\xi_{k'}
      + \Bigl(\frac{\alpha_{\c k}}{2Am} + \OO(n^{-1})
         \Bigr)\nu + \OO(n^{-3/2}),
         \displaybreak[0]\\
\mu &= \nu + \jsum
   \Bigl(\frac{\alpha_{j\c}}{2Amn} + \OO(n^{-2})\Bigr)\zeta_j
     + \ksum \Bigl( \frac{\alpha_{\c k}}{2Amn} +\OO(n^{-2})\Bigr)\xi_k
        + \OO(n^{-1})\nu,
\end{align*}
for $1\leq j\leq m-1$, $1\leq k\leq n-1$.

\medskip

The transformation $T_3^{-1}$ perturbs the region of integration in
an irregular  fashion that we must bound.  From the explicit form of
$T_3$ above, we have
\begin{align*}
\sigma_j &= \zeta_j + \sum_{j' = 1}^{m-1} \OO(n^{-3/2}) \zeta_{j'} +
               \ksum \OO(n^{-3/2}) \xi_k + \OO(n^{-1/2}) \nu+ \OO(n^{-3/2})
         = \zeta_j + \OO(n^{-1}),\\
\tau_k &= \xi_k + \sum_{j = 1}^{m-1} \OO(n^{-3/2}) \zeta_{j} +
                          \sum_{k'=1}^{n-1} \OO(n^{-3/2}) \xi_{k'} 
                          + \OO(n^{-1/2}) \nu+ \OO(n^{-3/2})
       = \xi_k + \OO(n^{-1})
\end{align*}
for $1\leq j\leq m-1$, $1\leq k\leq n-1$,
so $\sigmavec,\tauvec$ are only slightly different from $\zetavec,\xivec$.

For $\mu$ versus $\nu$ we have
\begin{align*}
  \mu &= \nu  + O(n^{-1+2\eps}/A) + O(m^{-1+2\eps}/A) \\
          &= \nu + o\( (mn)^{-1/2+2\eps}\),
\end{align*}
where the second step requires \eqref{oldassumptions}.
This shows that the bound $\abs{\nu}\le  (mn)^{-1/2+2\eps}$
is adequately covered by
$\abs\mu \le 2(mn)^{-1/2 + 2\eps}$.

For $1\leq h\leq 4$, define
\[ \mu_h = \jsum {\sigma_j}^h,\quad \nu_h = \ksum
            {\tau_k}^h.\]
{}From \eqref{zetaxibox}, we see that
$\abs{\pi_1} \le m^{1/2} n^{-1/2+ \eps}$ and
$\abs{\rho_1} \le m^{-1/2 + \eps}n^{1/2}$ are the remaining
constraints that define the region of integration.  We next apply
these constraints to bound $\mu_1$ and~$\nu_1$.
{}From the explicit form of $T_3$, we have
\begin{align}
 \mu_1 &= \pi_1
  + \jsum \sum_{j'=1}^{m-1}
  \Bigl(\frac{c(\alpha_{j\c}+
         \alpha_{j'\c})}{2An}+ \OO(n^{-2})\Bigr)\zeta_{j'}\notag 
    \\
  &  \kern10mm {}
   + \jsum\ksum
  \Bigl(\frac{\alpha_{jk} + d\alpha_{j\c} + c\alpha_{\c k}}
                      {2An} + \OO(n^{-2})\Bigr) \xi_k
  + \jsum \Bigl(\frac{\alpha_{j\c}}{2An}
  + \OO(n^{-1})\Bigr)\nu + \OO(n^{-1/2}) \notag \displaybreak[0]\\
 &= \pi_1 + \frac{c\alpha_{\c\c}}{2An} m^{1/2} n^{-1/2+\eps}
 + \frac{d\alpha_{\c\c}}{2An} m^{-1/2+\eps}n^{1/2}
      + \frac{\alpha_{\c\c}}{2An} \nu \notag \\
& \kern10mm  + \(1+c(m-1)\)\ksum\frac{\alpha_{\c k}}{2An}\xi_k
 + \frac{c(m-1)}{2An}\sum_{j'=1}^{m-1} \alpha_{j'\c}\zeta_{j'}
  + \OO(n^{-1/2}) \notag \displaybreak[0]\\
 &= \pi_1
  + \frac{c(m-1)}{2An} \sum_{j'=1}^{m-1} \alpha_{j'\c}\zeta_{j'}
  + \OO(n^{-1/2}) \label{mu1pi1}\\
  &= \pi_1 +O(A^{-1}mn^{-1+2\eps})\notag  \\[0.5ex]
  &= \pi_1 + o(m^{1/2}n^{-1/2+5\eps/2}).\notag
\end{align}
To derive the above we have used $1+c(m-1)=m^{1/2}$ and the bounds we
have established on the various variables.  For the last step, we
need \eqref{oldassumptions},
which implies that
$A^{-1}mn^{-1+2\eps}=o(m^{1/2}n^{-1/2+5\eps/2})$.

Since our region of integration has
$\abs{\pi_1}\le m^{1/2}n^{-1/2+\eps}$, we see that this implies
the bound $\abs{\mu_1}\le m^{1/2}n^{-1/2+3\eps}$.
By a parallel argument, we have
\[\nu_1 = \rho_1 + o(m^{-1/2+5\eps/2}n^{1/2}),\]
which implies
$\abs{\nu_1}\le n^{1/2}m^{-1/2+3\eps}$.
Putting together all the bounds we have derived, we see that
\[T_3^{-1}(T_2^{-1}(\S)) \subseteq \Q\cap\M,\]
where
\begin{align*}
  \Q &= \{\, \abs{\sigma_j}\le 2n^{-1/2 + \eps}, j=1,\ldots, m-1\,\}
   \cap \{\, \abs{\tau_k}\le 2m^{-1/2 + \eps}, k=1,\ldots, n-1\,\}\\
  & \hspace*{1cm} \cap \{ \abs\mu \le 2(mn)^{-1/2+2\eps}\,\},\\
  \M &= \{\, \abs{\mu_1}\le m^{1/2} n^{-1/2+3\eps}\,\} \cap
        \{\, \abs{\nu_1}\le n^{1/2} m^{-1/2+3\eps}\}.
\end{align*}

Now define
\begin{align}\label{RSprime}
\begin{split}
 \S' &=T_2(T_3(\Q\cap\M)), \\
 \R' &= T_1^*(\S').
\end{split}
\end{align}
We have proved that $\S'\supseteq\S$.
Also notice that $\R'$ is contained in the region
defined by the inequalities~\eqref{bigR}.
As we forecast at that time, our estimates of the integrand
have been valid inside this expanded region.  It remains
to apply the transformation $T_3^{-1}$ to the
integrand~\eqref{Helper} so that we have it in terms of
$(\sigmavec,\tauvec,\mu)$.  The explicit form of $T_3^{-1}$
is similar to the explicit form for $T_3$, namely:
\begin{align*}
\zeta_j &= \sigma_j - \sum_{j'=1}^{m-1} \Bigl(\frac{c(\alpha_{j\c} +
        \alpha_{j'\c})}{2An} + \OO(n^{-2})\Bigr) \sigma_{j'} -
        \ksum \Bigl(\frac{\alpha_{jk} + d\alpha_{j\c} +
          c\alpha_{\c k}}{2An} + \OO(n^{-2})\Bigr) \tau_k \\
      & \hspace*{2cm} {} -
         \Bigl(\frac{\alpha_{j\c}}{2An} + \OO(n^{-1})\Bigr)\mu
          +\OO(n^{-3/2}),
       \displaybreak[0]\\
\xi_k &= \tau_k - \jsum \Bigl(\frac{\alpha_{jk} -
      d\alpha_{j\c} + c\alpha_{\c k}}{2Am} + \OO(n^{-2})\Bigr)\sigma_j
    - \sum_{k'=1}^{n-1} \Bigl(\frac{d(\alpha_{\c k} + \alpha_{\c k'})}
  {2Am} + \OO(n^{-2})\Bigr)\tau_{k'} \\
      & \hspace*{2cm} {} - \Bigl(\frac{\alpha_{\c k}}{2Am} + \OO(n^{-1})
         \Bigr)\mu +\OO(n^{-3/2}),\\
\nu &= \mu - \jsum
   \Bigl(\frac{\alpha_{j\c}}{2Amn} + \OO(n^{-2})\Bigr)\sigma_j
     - \ksum \Bigl( \frac{\alpha_{\c k}}{2Amn} +\OO(n^{-2})\Bigr)\tau_k
        + \OO(n^{-1})\mu,
\end{align*}
for $1\leq j\leq m-1$, $1\leq k\leq n-1$.
In addition to the relationships between the old and new variables
that we proved before, we can note that $\pi_2=\mu_2+\OO(n^{-1/2})$,
$\rho_2=\nu_2+\OO(n^{-1/2})$, $\pi_3=\mu_3+\OO(n^{-1})$,
$\rho_3=\nu_3+\OO(n^{-1})$,  $\pi_4=\mu_4+\OO(n^{-3/2})$,
and $\rho_4=\nu_4+\OO(n^{-3/2})$.

\medskip
The quadratic part of $L_1$, which we called $Q$ in~\eqref{Qvalue},
loses its off-diagonal parts according to our design of $T_3$.
Thus, what remains is
\begin{align*}
 -Amn\mu^2 &- \jsum \(An+(1+2c)\alpha_{j\c}-Ax'_j\)\sigma_j^2
   - \ksum \(Am+(1+2d)\alpha_{\c k}-Ay'_k\)\tau_k^2\\
     &= -Amn\mu^2 - An\mu_2 - Am\nu_2 \\
     &{~\quad}
         - \jsum (\alpha_{j\c}-Ax'_j) \sigma_j^2
         - \ksum (\alpha_{\c k}-Ay'_k) \tau_k^2 + \OO(n^{-1/2}).
\end{align*}

Next consider the cubic terms of $L_1$.  These are
\begin{align*}
   &{}-3i A_3 n\nu \pi_2 - 3iA_3 m\nu \rho_2
   - i A_3 n\pi_3 - i A_3 n\rho_3 \\
&\qquad{}- 3i A_3 cn \pi_1\pi_2 - 3iA_3 d n \rho_1\rho_2 
- i \jsum\ksum \beta_{jk}(\zeta_j+\xi_k)^3.
\end{align*}
We calculate the following in $\Q\cap\M$:
\begin{align}
 -3i A_3 n\nu \pi_2 &= -3i A_3 n\mu\mu_2
       + \frac{3iA_3\mu_2}{2Am}\biggl(\,\jsum \alpha_{j\c}\sigma_j
                    + \ksum \alpha_{\c k}\tau_k\biggr) + \OO(n^{-1/2}), \notag\\
  -i A_3 n\pi_3 &= -i A_3 n\mu_3 + \frac{3iA_3}{2A}\biggl(\,
       \sum_{j,j'=1}^{m-1}
         c(\alpha_{j\c}+\alpha_{j'\c})\sigma_j^2\sigma_{j'}, \notag \\
    &\kern40mm{}   + \jsum\ksum
           (\alpha_{jk}+d\alpha_{j\c}+c\alpha_{\c k})
             \sigma_j^2\tau_k\biggr) + \OO(n^{-1/2}), \notag 
             \displaybreak[0]\\
  - 3i A_3 cn \pi_1\pi_2 &= - 3i A_3 cn \mu_1\mu_2
     + \frac{3iA_3c^2m\mu_2}{2A}
           \jsum \alpha_{j\c}\sigma_j + \OO(n^{-1/2}) \label{pi1pi2}, \\
  - i \jsum\ksum \beta_{jk}&(\zeta_j+\xi_k)^3 =
     - i \jsum\ksum \beta_{jk}(\sigma_j+\tau_k)^3
         + \OO(n^{-1/2}),\notag
\end{align}
and the remaining cubic terms are each parallel to one of those.
The proof of \eqref{pi1pi2} is similar to the proof of \eqref{mu1pi1}.

Finally we come to the quartic part of $L_1$, which is
\[6A_4\pi_2\rho_2 + A_4n\pi_4 + A_4m\rho_4
  = 6A_4\mu_2\nu_2 + A_4n\mu_4 + A_4m\nu_4+ \OO(n^{-1/2}).\]

In summary, the value of the integrand for
$(\sigmavec,\tauvec,\mu)\in\Q\cap\M$ is
$\exp\(L_2+\OO(n^{-1/2})\)$, where
\begin{align*}
\begin{split}
   L_2 &= -Amn\mu^2 - An\mu_2 - Am\nu_2
         - \jsum (\alpha_{j\c}-Ax'_j) \sigma_j^2
         - \ksum (\alpha_{\c k}-Ay'_k) \tau_k^2
         + 6A_4\mu_2\nu_2 \\
         &\quad{}+ A_4n\mu_4 + A_4m\nu_4
      - i A_3 n\mu_3 - i A_3 m\nu_3
               - 3i A_3 cn \mu_1\mu_2 - 3i A_3 dm \nu_1\nu_2\\
       &\quad{} - 3 i A_3 n \mu\mu_2 - 3 i A_3 m \mu\nu_2
                - i \jsum \beta_{j\c}\sigma_j^3
               - i \ksum \beta_{\c k}\tau_k^3 \\
       &\quad{} + i\sum_{j,j'=1}^{m-1} g_{jj'}\sigma_j\sigma_{j'}^2
             + i\sum_{k,k'=1}^{n-1} h_{kk'}\tau_k\tau_{k'}^2
             + i\jsum\ksum
                \(u_{jk}\sigma_j\tau_k^2 + v_{jk}\sigma_j^2\tau_k\),
\end{split}
\end{align*}
with
\begin{align*}
    g_{jj'} &= \frac{3A_3}{2Am}
                    \((1+cm+c^2m^2)\alpha_{j\c}+cm\alpha_{j'\c}\)
                    = O(n^{-1/2+\eps}), \\[0.4ex]
    h_{kk'} &= \frac{3A_3}{2An}
                    \((1+dn+d^2n^2)\alpha_{\c k}+dn\alpha_{\c k'}\)
                    = O(m^{-1/2+\eps}), \displaybreak[1]\\[0.4ex]
    u_{jk} &= \frac{3A_3}{2An}
                \( n\alpha_{jk} + (1+dn)\alpha_{j\c} + cn\alpha_{\c k}\)
                   - 3\beta_{jk} = O(m^{-1/2+2\eps}+n^{-1/2+2\eps}), \\[0.4ex]
    v_{jk} &= \frac{3A_3}{2Am}
                \( m\alpha_{jk} + (1+cm)\alpha_{\c k} + dm\alpha_{j\c}\)
                           - 3\beta_{jk} = O(m^{-1/2+2\eps}+n^{-1/2+2\eps}) .
\end{align*}
Note that the $O(\cdot)$ estimates in the last four lines are uniform over
$j,j',k,k'$.

\nicebreak
\subsection{Estimating the main part of the integral}\label{s:complete}

Define $E_2=\exp(L_2)$.  We have shown that the value of the
integrand in
$\Q\cap\M$ is $E_1=E_2\(1+\OO(n^{-1/2})\)$.
Denote the complement of the region $\M$ by $\M^c$.   We can
approximate our integral as follows:
\begin{align}
\int_{\Q \cap \M} E_1
  &= \int_{\Q \cap \M} E_2 + \OO(n^{-1/2})\int_{\Q \cap \M} \abs{E_2} \notag\\[0.4ex]
  &= \int_{\Q \cap \M} E_2 + \OO(n^{-1/2})\int_{\Q} \,\abs{E_2} \notag\\[0.4ex]
  &= \int_{\Q} E_2
        + O(1)\int_{\Q\cap \M^c} \abs{E_2} + \OO(n^{-1/2})\int_{\Q} \,\abs{E_2}.
                  \label{recipe}
\end{align}
It suffices to estimate the value of each integral in \eqref{recipe}.
This can be done using the same calculation as in Section~4.3
of~\cite{CGM}, using $\hat\alpha_{jk}=\alpha_{jk}-A\h_{jk}$
in place of the variable~$\alpha_{jk}$ used in that paper.
A potential problem with this analogy is that the variable
$\alpha_{jk}$ used in~\cite{CGM}
has the property $\alpha_{jk}=\OO(n^{-1/2})$,
whereas it is not true that $\hat\alpha_{jk}=\OO(n^{-1/2})$.
However, a careful look at Section~4.3 of~\cite{CGM} confirms that
only the properties $\hat\alpha_{j\c}=\alpha_{j\c}-Ax'_j=\OO(n^{1/2})$,
$\hat\alpha_{\c k}=\alpha_{\c k}-Ay'_k=\OO(n^{1/2})$, and the bounds
on $g_{jj'}, h_{kk'}, u_{jk}, v_{jk}$,
are required.

The result is that
\begin{equation}\label{QE2a}
\begin{split}
 \int_\Q E_2
  &= \Bigl(\frac{\pi}{Amn}\Bigr)^{1/2}
        \Bigl(\frac{\pi}{An}\Bigr)^{\!(m-1)/2}
        \Bigl(\frac{\pi}{Am}\Bigr)^{\!(n-1)/2} \\
  &\qquad{}\times\exp\biggl(  - \frac{9A_3^2}{8A^3} + \frac{3A_4}{2A^2}
         +\Bigl(\frac mn+\frac nm\Bigr)
         \Bigl( \frac{3A_4}{4A^2}-\frac{15A_3^2}{16A^3}\Bigr) \\
   &\kern23mm{} 
      - \Bigl(\frac{1}{2Am}+\frac{1}{2An}\Bigr)\hat\alpha_{\c\c}
      + \frac{1}{4A^2m^2}\ksum (\hat\alpha_{\c k})^2 \\
   &\kern23mm{}
      + \frac{1}{4A^2n^2}\jsum (\hat\alpha_{j\c})^2
      + \OO(n^{-b}) \biggr),
\end{split}
\end{equation}
where $b$ is specified in Theorem~\ref{bigtheorem}. 

Using \eqref{rad10} and the conditions of
Theorem~\ref{bigtheorem}, we calculate that
\begin{align*}
\hat\alpha_{\c\c} &= 
    -\frac12\Bigl(\frac{R_2}{n}+\frac{C_2}{m}\Bigr)
    -\dfrac12\lambda^2 X +\OO(n^{1/2}), \displaybreak[0]\\
\jsum (\hat\alpha_{j\c})^2 
  &= \dfrac14 (1-2\lambda)^2R_2+\OO(n^{3/2}), \\
\ksum (\hat\alpha_{\c k})^2 
  &= \dfrac14 (1-2\lambda)^2 C_2+\OO(n^{3/2}).
\end{align*}
Substituting these values into \eqref{QE2a} together with the
actual values of $A, A_3, A_4$, we conclude that
\begin{equation}\label{QE2}
\begin{split}
 \int_\Q E_2
  &= \Bigl(\frac{\pi}{Amn}\Bigr)^{1/2}
        \Bigl(\frac{\pi}{An}\Bigr)^{\!(m-1)/2}
        \Bigl(\frac{\pi}{Am}\Bigr)^{\!(n-1)/2} \\
   &{\kern7mm}\times \exp\biggl(
         - \frac12
         - \frac{1-2A}{24A}\Bigl(\frac mn+\frac nm\Bigr)
         + \frac{1-4A}{16A^2}\Bigl(\frac{R_2}{n^2}+\frac{C_2}{m^2}\Bigr)  
   \\
      &{\kern22mm}
      + \frac{R_2+C_2}{4Amn}
      + \frac{\lambda^2 X}{4A}\Bigl(\frac{1}{m}+\frac{1}{n}\Bigr)
       + O(n^{-b})  \biggr).
\end{split}
\end{equation}

By the same argument as in~\cite{CGM}, the other two terms
in~\eqref{recipe} have value $O(n^{-b})  \int_Q E_2$.
Multiplying \eqref{QE2}
by the Jacobians of the transformations~$T_2$ and~$T_3$,
we have proved the following.

\begin{lemma}\label{Jintegral}
The region $\S'$ given by~\eqref{RSprime} contains $\S$ and
\begin{align*}
 \int_{\S'} G(\thetahatvec,\phihatvec,\nu)
     \, d\thetahatvec  d\phihatvec  d\nu
     &= (mn)^{-1/2} \Bigl(\frac{\pi}{Amn }\Bigr)^{\!1/2}
                    \Bigl(\frac{\pi}{An  }\Bigr)^{(m-1)/2}
                    \Bigl(\frac{\pi}{Am  }\Bigr)^{(n-1)/2} \\
     &{\kern7mm}\times \exp\biggl(
         - \frac12
         - \frac{1-2A}{24A}\Bigl(\frac mn+\frac nm\Bigr)
         + \frac{1-4A}{16A^2}\Bigl(\frac{R_2}{n^2}+\frac{C_2}{m^2}\Bigr)  
   \\
      &{\kern20mm}
      + \frac{R_2+C_2}{4Amn}
      + \frac{\lambda^2 X}{4A}\Bigl(\frac{1}{m}+\frac{1}{n}\Bigr)
       + O(n^{-b})  \biggr).
\end{align*}
\end{lemma}

\nicebreak
\subsection{Bounding the remainder of the integral}\label{s:boxing}

In the previous subsection, we estimated the value of the integral
$I_{\R'}(\svec,\tvec,\XX)$, which is the same as $I(\svec,\tvec,\XX)$
except that it is restricted to a certain region $\R'\supseteq\R$.
In this subsection, we
extend this to an estimate of $I(\svec,\tvec,\XX)$ by showing that
the remainder of the region of integration contributes
negligibly.

For $1\leq j\leq m$, $1\leq k\leq n$, let
$A_{jk} = A + \alpha_{jk}=\tfrac12\lambda_{jk}(1-\lambda_{jk})$
(recall \eqref{AlBetGamDef}),
and define $\Amin=\min_{jk} A_{jk} = A+\OO(n^{-1/2})$.
We begin with two technical lemmas whose proofs are omitted,
and a well-known bound of Hoeffding.
\begin{lemma}\label{fbnd}
\[
 \abs{F(\thetavec,\phivec)} = \prod_{jk\in\barXX} f_{jk}(\t_j + \p_k),
\]
where
\[f_{jk}(z) =\ssqrt{1-4A_{jk}(1-\cos z)}\,.\]
Moreover, for all real $z$,
\[
 0\le f_{jk}(z) \le \exp\( -A_{jk} z^2 + \dfrac1{12} A_{jk} z^4\).
    \quad\qedsymbol
\]
\end{lemma}

\begin{lemma}\label{ibnd}
For all\/ $c>0$,
\[
 \int_{-8\pi/75}^{8\pi/75} \exp\( c (-x^2 + \dfrac73 x^4)\)\,dx
 \le \sqrt{\pi/c}\, \exp(3/c).   \quad\qedsymbol
\]
\end{lemma}

\begin{lemma}[{\cite{hoeffding}}]\label{hoeffding}
Let $X_1,\ldots,X_N$ be independent random
variables such that $\expect X_i=0$ and $\abs{X_i}\le M$ for all~$i$.
Then, for any $t\ge 0$,
\[
 \Prob\Bigl(\,\sum_{i=1}^N X_i \ge t\Bigr)
    \le \exp\biggl(-\frac{t^2}{2NM^2}\biggr).\quad\qedsymbol
\]
\end{lemma}

\begin{lemma}\label{boxing}
  Let $F(\thetavec,\phivec)$ be the integrand of $I(\svec,\tvec,\XX)$
  as defined in~\eqref{start}.
  Then, under the conditions of Theorem~\ref{bigtheorem},
  \[
   \int_{\R^c} \abs{F(\thetavec,\phivec)}\,d\thetavec d\phivec
    = O(n^{-1}) \int_{\R'} F(\thetavec,\phivec)\,d\thetavec d\phivec,
  \]
where $\R^c$ denotes the complement of $\R$.
\end{lemma}

\begin{proof}
Our approach will be to bound $\int \,\abs{F(\thetavec,\phivec)}$ over
a variety of regions whose union covers $\R^c$.
To make the comparison of these bounds with
$\int_{\R'} F(\thetavec,\phivec)$ easier, we note that
\begin{equation}\label{I0I1}
   \int_{\R'} F(\thetavec,\phivec)\,d\thetavec d\phivec
     = \exp\(A^{-1}O(m^\eps+n^\eps)\) I_0 
     = \exp\(O(m^{3\eps}+n^{3\eps})\) I_1,
\end{equation}
where
\begin{align*}
   I_0 &= \Bigl(\frac{\pi}{A_{\b\b}}\Bigr)^{1/2}\,
   \prod_{j=1}^m   \Bigl(\frac{\pi}{A_{j\b} }\Bigr)^{\!1/2}\,
          \prod_{k=1}^n   \Bigl(\frac{\pi}{A_{\b k} }\Bigr)^{\!1/2}
          \!, \\
   I_1 &= \Bigl(\frac{\pi}{An}\Bigr)^{\!m/2} \Bigl(\frac{\pi}{Am}\Bigr)^{\!n/2}\!.
\end{align*}
To see this, expand 
\[
 A_{j\b} = An+\alpha_{j\b}=An \exp\biggl( \frac{\alpha_{j\b}}{An}
   - \frac{\alpha_{j\b}^2}{2A^2n^2} +\cdots\,\biggr),\]
and similarly for $A_{\b k}$, and compare the result
to Lemma~\ref{Jintegral} using the assumptions of
Theorem~\ref{bigtheorem}.  It may help to recall the calculation
following~\eqref{QE2a}.

\smallskip

Take $\kappa=\pi/300$ and define $w_0,w_1,\ldots,w_{299}$ by
$w_\ell = 2\ell\kappa$.
For any $\ell$,
let $\S_1(\ell)$ be the set of
$(\thetavec,\phivec)$ such that $\t_j\in [w_\ell-\kappa,w_\ell+\kappa]$ for
at least $\kappa m/\pi$ values of~$j$ and
$\p_k\notin [-w_\ell-2\kappa,-w_\ell+2\kappa]$ for at
least $n^\eps$ values of~$k$.
For $(\thetavec,\phivec)\in\S_1(\ell)$, $\theta_j+\phi_k\notin[-\kappa,\kappa]$
for at least $\kappa (m-O(m^\eps))n^\eps/\pi$ pairs $(j,k)$
with $\h_{jk}=0$ so, by Lemma~\ref{fbnd},
$\abs{F(\thetavec,\phivec)}\le \exp(-c_1\Amin mn^\eps)$
for some $c_1>0$ which is independent of~$\ell$.

Next  define $\S_2(\ell)$ to be the set of
$(\thetavec,\phivec)$ such that
$\t_j\in [w_\ell-\kappa,w_\ell+\kappa]$ for
at least $\kappa m/\pi$ values of~$j$,
$\p_k\in[-w_\ell-2\kappa,-w_\ell+2\kappa]$ for at least $n-n^\eps$ values of~$k$
and
$\t_j\notin [w_\ell-3\kappa,w_\ell+3\kappa]$ for at least $m^\eps$ values of~$j$.
By the same argument with the roles of $\thetavec$
and $\phivec$ reversed,
$\abs{F(\thetavec,\phivec)}\le \exp(-c_2\Amin m^\eps n)$
for some $c_2>0$ independent of~$\ell$ when
$(\thetavec,\phivec)\in\S_2(\ell)$.

Now define $\R_1(\ell)$ to be the set of
pairs $(\thetavec,\phivec)$ such that
$\t_j\in [w_\ell-3\kappa,w_\ell+3\kappa]$ for at least $m-m^\eps$
values of~$j$, and $\p_k\in [-w_\ell-3\kappa,-w_\ell+3\kappa]$ for
at least $n-n^\eps$ values of~$k$.
By the pigeonhole principle, for any $\thetavec$ there is some~$\ell$
such that $[w_\ell-\kappa,w_\ell+\kappa]$ contains at least $\kappa m/\pi$
values of $\t_j$.  Therefore,
\[
    \Bigl(\,\bigcup_{\ell=0}^{299} \, \R_1(\ell)\!\Bigr)^c 
        \subseteq \bigcup_{\ell=0}^{299}  \,\(\S_1(\ell)\cup\S_2(\ell)\).
\]
Since the total volume of $\(\,\bigcup_\ell\R_1(\ell)\)^c$ is at most $(2m)^{m+n}$,
we find that for some $c_3>0$,
\begin{align}
 \int_{(\bigcup_\ell \R_1(\ell))^c} \,&
    \abs{F(\thetavec,\phivec)}\,d\thetavec d\phivec \notag\\
  & \le
     (2\pi)^{m+n} \( \exp( -c_3 \Amin m n^\eps)
         +\exp( -c_3 \Amin m^\eps n)\) \notag\\
  & \le e^{-n} I_1.\label{R1C}
\end{align}

We are left with $(\thetavec,\phivec)\in\bigcup_\ell\R_1(\ell)$.
If we subtract $w_\ell$ from each $\t_j$ and add $w_\ell$ to each $\p_k$ the
integrand $F(\thetavec,\phivec)$ is unchanged, so we can assume
for convenience that $\ell=0$ and that
$(\thetavec,\phivec)\in\R_1=\R_1(0)$.
The bounds we obtain on parts of the integral we seek to reject will
be at least $1/300$ of the total and thus be of the right order of
magnitude.  We will not mention this point again.

For a given $\thetavec$, partition $\{1,2,\ldots, m\}$ into sets
$J_0=J_0(\thetavec)$, $J_1 = J_1(\thetavec)$ and $J_2 = J_2(\thetavec)$,
containing the indices $j$ such that
$\abs{\t_j} \leq 3\kappa$, $3\kappa  < \abs{\t_j}\leq 15\kappa$
and $\abs{\t_j} > 15\kappa$, respectively.
Similarly partition $\{1,2,\ldots,n\}$ into
$K_0=K_0(\phivec)$, $K_1=K_1(\phivec)$ and $K_2=K_2(\phivec)$.
The value of $\abs{F(\thetavec,\phivec)}$ can now be bounded using
\begin{align*}
f_{jk}(\t_j+\p_k) & \\
 &\kern-14mm{}\le \begin{cases}
    \,\exp\(-\Amin(\t_j+\p_k)^2+\tfrac1{12}\Amin(\t_j+\p_k)^4\)
        & \textrm{if }(j,k) \in (J_0\cup J_1)\times(K_0\cup K_1), \\[0.8ex]
    \,\ssqrt{\vrule width0pt depth0.3ex
          1-4\Amin(1-\cos(12\kappa))} \le e^{-\Amin/64}
        & \textrm{if }(j,k) \in (J_0\times K_2)\cup (J_2\times K_0), \\[0.3ex]
    \,1   & \textrm{otherwise}.
 \end{cases}
\end{align*}
Let $I_2(m_2,n_2)$ be the contribution to
$\int_{\R_1} \abs{F(\thetavec,\phivec)}$
of those $(\thetavec,\phivec)$ with
$\card{J_2}=m_2$ and $\card{K_2}=n_2$.
Recall that $\card{J_0}> m-m^\eps$ and $\card{K_0}> n-n^\eps$.
We have
\begin{align}\label{I2bnd}
\begin{split}
  I_2(m_2,n_2) 
 &\le \binom{m}{m_2}\binom{n}{n_2}(2\pi)^{m_2+n_2} \\
       &\kern7mm{}\times \exp\( -\dfrac1{64} \Amin(n-O(n^\eps))m_2
               -\dfrac1{64} \Amin(m-O(m^\eps))n_2\) I'_2(m_2,n_2),
\end{split}
\end{align}
where
\[
  I'_2(m_2,n_2) =
   \int_{-15\kappa}^{15\kappa}\!\!\cdots \int_{-15\kappa}^{15\kappa}
    \!\!\!\exp\Bigl(- \Amin\sumppd_{jk\in\barXX}(\t_j+\p_k)^2 +
                \dfrac1{12}\Amin\sumppd_{jk\in\barXX}(\t_j+\p_k)^4\Bigr)
      \, d\thetavec' d\phivec',
\]
and the primes denote restriction to $j\in J_0\cup J_1$
and $k\in K_0\cup K_1$, in the case of the summations in addition
to the restriction given by the summation limits.
Write $m'=m-m_2$ and $n'=n-n_2$ and define
$\avtheta'=(m')^{-1}\sumpp_j\t_j$,
$\breve\theta_j=\t_j-\avtheta'$ for $j\in J_0\cup J_1$,
$\avphi'=(n')^{-1}\sumpp_k\p_k$,
$\breve\phi_k=\p_k-\avphi'$ for $k\in K_0\cup K_1$,
$\nu'=\avphi'+\avtheta'$ and
$\psi'=\avtheta'-\avphi'$.
Change variables from $(\thetavec',\phivec')$ to
$\{\breve\theta_j \mathrel| j\in J_3\} 
\cup \{\breve\phi_k\mathrel| k\in K_3\}
\cup \{\nu',\psi'\}$, where $J_3$ is some subset of $m'-1$ elements
of $J_0\cup J_1$ and $K_3$ is some subset of $n'-1$ elements
of $K_0\cup K_1$.  {}From Subsection~\ref{s:integral} we know that the
Jacobian of this transformation is $m'n'/2$.
The integrand of $I'_2$ can now be bounded using
\[ \sumppd_{jk\in\barXX}(\t_j+\p_k)^2 =
   (n'-O(n^\eps))\sumppd_j\breve\theta_j^2 
   + (m'-O(m^\eps))\sumppd_k\breve\phi_k^2
   + (m'n'-O(X))\nu'{}^2
\]
and
\[ \sumppd_{jk\in\barXX}(\t_j+\p_k)^4 \le 
     27n'\sumppd_j\breve\theta_j^4 + 27m'\sumppd_k\breve\phi_k^4
     + 27m'n'\nu'{}^4.
\]
The latter follows from the inequality
$(x+y+z)^4 \le 27(x^4+y^4+z^4)$ valid for all $x,y,z$.
Therefore,
\begin{align*}
I'_2(m_2,n_2) \le \frac{O(1)}{m'n'} \int_{-30\kappa}^{30\kappa}
    \int_{-30\kappa}^{30\kappa}\!\!\!\cdots\!\int_{-30\kappa}^{30\kappa}
     \!\exp\Bigl(& \Amin (n'-O(n^\eps)) \sumppd_j g(\breve\theta_j)\\[-1ex]
    &{\!}+ \Amin (m'-O(m^\eps))\sumppd_k g(\breve\phi_k) \\
    &{\!}+ \Amin (m'n'-O(X)) g(\nu') \Bigr)\,
    d\breve\theta_{j\in J_3}d\breve\phi_{k\in K_3}\, d\nu',
\end{align*}
where $g(z) = -z^2+\tfrac73 z^4$.  Since $g(z)\le 0$ for
$\abs{z}\le 30\kappa$, and we only need
an upper bound, we can restrict the summations in the
integrand to $j\in J_3$ and
$k\in K_3$.  The integral now separates into $m'+n'-1$
one-dimensional integrals
and Lemma~\ref{ibnd} (by monotonicity) gives that
\begin{align*}
I'_2(m_2,n_2) &= O(1)
     \frac{\pi^{(m'+n')/2}}
     { \Amin^{(m'+n'-1)/2}
         (m'-O(m^\eps))^{n'/2-1} (n'-O(n^\eps))^{m'/2-1}} \\[0.5ex]
    &{\qquad}\times \exp\( O(m'/(\Amin n') + n'/(\Amin m'))\).
\end{align*}
Applying \eqref{I0I1} and \eqref{I2bnd}, we find that
\begin{equation}\label{I2mn}
 \mathop{\sum_{m_2=0}^{m^\eps} \;
                \sum_{n_2=0}^{n^\eps}}\displaylimits_{m_2+n_2\ge 1}
 I_2(m_2,n_2) = O\(e^{-c_4A m}+e^{-c_4A n}\) I_1
\end{equation}
for some $c_4>0$.

\medskip

We have now bounded contributions to the integral of
$\abs{F(\thetavec,\phivec)}$ from everywhere outside
the union of 300 equivalent translates of $\X-\R$, where
\[
\X = \bigl\{\, (\thetavec,\phivec) \bigm|
  \abs{\t_j},\abs{\p_k} \le15\kappa~\text{for}~1\le j\le m, 1\le k\le n\, \bigr\}.
\]
By Lemma~\ref{fbnd}, we have for $(\thetavec,\phivec)\in\Chat^{m+n}$
(which includes~$\X$) that
\[
\abs{F(\thetavec,\phivec)} \le\exp\Bigl(
 - \sum_{jk\in\barXX} A_{jk} (\that_j+\phat_k+\nu)^2
  + \dfrac1{12}\sum_{j=1}^m\sum_{k=1}^n A_{jk} (\that_j+\phat_k+\nu)^4\Bigr),
\]
where $\that_j = \theta_j - \avtheta$, $\phat_k = \phi_k - \avphi$
and $\nu=\avtheta+\avphi$.  As before, the integrand is independent
of $\psi=\avtheta-\avphi$ and our notation will tend to ignore $\psi$
for that reason; for our bounds it will suffice to remember that
$\psi$ has a bounded range.

We proceed by exactly diagonalizing the $(m+n+1)$-dimensional
quadratic form.  Since $\sum_{j=1}^m\that_j = \sum_{k=1}^n\phat_k=0$,
we have
\begin{align*}
 \sum_{jk\in\barXX}
   A_{jk} (\that_j+\phat_k+\nu)^2
  &= \sum_{j=1}^m A_{j\b|\barXX\,}\that_j^2
    + \sum_{k=1}^n A_{\b k|\barXX\,}\phat_k^2 + A_{\b\b|\barXX\,}\nu^2\\
    &\quad{}+ 2\sum_{j=1}^m\sum_{k=1}^n
      (\alpha_{jk}-A_{jk}\h_{jk})\that_j\phat_k \\
       &\quad{} + 2\nu\sum_{j=1}^m (\alpha_{j\b}-A_{j\b|\XX})\that_j 
        + 2\nu\sum_{k=1}^n (\alpha_{\b k}-A_{\b k|\XX})\phat_k.
\end{align*}
This is almost diagonal, because $\alpha_{jk} = \OO(n^{-1/2})$,
$A_{j\b|\XX}=\OO(1)$, $A_{\b k|\XX}=\OO(1)$.  The coefficients
$-2A_{jk}\h_{jk}$ can be larger but only in the $\OO(n)$ places
where $\h_{jk}=1$.
We can make the quadratic form exactly diagonal using
the slight additional transformation $(\II+\UU^{-1}\YY)^{-1/2}$
described by Lemma~\ref{diagonalize}, where $\UU$ is a diagonal
matrix with diagonal entries $A_{j\b|\barXX}$,
$A_{\b k|\barXX}$ and~$A_{\b\b|\barXX}$.
The matrix $\YY$ has zero diagonal and other entries of
magnitude $\OO(n^{-1/2})$ apart from the row and column
indexed by $\nu$, which
have entries of magnitude~$\OO(n^{1/2})$,
and the $\OO(n)$ just-mentioned entries of order $\OO(1)$.
By the same argument as used in Subsection~\ref{s:diagonalize},
all eigenvalues of $\UU^{-1}\YY$ have  magnitude~$\OO(n^{-1/2})$, so the
transformation is well-defined.
The new variables
$\{\hat\vartheta_j\}$, $\{\hat\varphi_k\}$ and $\varnu$
are related to the old by
\[
(\that_1,\ldots,\that_m,\phat_1,\ldots,\phat_n,\nu)^T
 =(\II+\UU^{-1}\YY)^{-1/2}(\hat\vartheta_1,\ldots,\hat\vartheta_m,
 \hat\varphi_1,\ldots,\hat\varphi_n,\varnu)^T.
\]
We will keep the variable~$\psi$ as a variable of
integration but, as noted before, our notation will generally ignore~it.

More explicitly, for some $d_1,\ldots,d_m,d'_1,\ldots,d'_n=\OO(n^{-3/2})$,
we have uniformly over\\ $j=1,\ldots,m$, $k=1,\ldots,n$ that
\begin{align}\label{varrels}
\begin{split}
  \that_j &= \hat\vartheta_j + \sum_{q=1}^m \OO(n^{-2})\hat\vartheta_q
           + \sum_{k=1}^n \OO(n^{-3/2}+n^{-1}\h_{jk})\hat\varphi_k 
           + \OO(n^{-1/2})\varnu, \\
  \phat_k &= \hat\varphi_k 
           + \sum_{j=1}^m \OO(n^{-3/2}+n^{-1}\h_{jk})\hat\vartheta_j
           + \sum_{q=1}^n \OO(n^{-2})\hat\varphi_q 
           + \OO(n^{-1/2})\varnu, \\
  \nu &= \varnu + \sum_{j=1}^m d_j\hat\vartheta_j
                         + \sum_{k=1}^n d'_k \hat\varphi_k + \OO(n^{-1})\varnu.
\end{split}
\end{align}
Note that the expressions $\OO(\cdot)$ in \eqref{varrels} represent values
that depend on $m,n,\svec,\tvec$ but not on $\{\hat\vartheta_j\},
\{\hat\varphi_k\},\varnu$.

The region of integration $\X$ is
$(m{+}n)$-dimensional.  In place of the variables
$(\thetavec,\phivec)$ we can use $(\thetahatvec,\phihatvec,\nu,\psi)$
by applying the identities $\that_m=-\sum_{j=1}^{m-1}\that_j$ and
$\phat_n=-\sum_{k=1}^{n-1}\phat_k$.
(Recall that $\thetahatvec$ and $\phihatvec$ don't include
$\that_m$ and $\phat_n$.)
The additional
transformation \eqref{varrels} maps the two just-mentioned identities
into identities that define $\hat\vartheta_m$ and $\hat\varphi_n$
in terms of $(\varthetahatvec,\varphihatvec,\varnu)$, where
$\varthetahatvec=(\hat\vartheta_1,\ldots,\hat\vartheta_{m-1})$
and $\varphihatvec=(\hat\varphi_1,\ldots,\hat\varphi_{n-1})$.
These have the form
\begin{equation}\label{mnvars}
\begin{split}
 \hat\vartheta_m &= -\sum_{j=1}^{m-1}\(1+\OO(n^{-1})\)\hat\vartheta_j
      + \sum_{k=1}^{n-1}\OO(n^{-1/2})\hat\varphi_k + \OO(n^{1/2})\varnu,\\
 \hat\varphi_n &= \sum_{j=1}^{m-1}\OO(n^{-1/2})\hat\vartheta_j 
       - \sum_{k=1}^{n-1}\(1+\OO(n^{-1})\)\hat\varphi_k + \OO(n^{1/2})\varnu.
\end{split}
\end{equation}
Therefore, we can now integrate over
$(\varthetahatvec,\varphihatvec,\varnu,\psi)$.
The Jacobian of the transformation from $(\thetavec,\phivec)$
to $(\thetahatvec,\phihatvec,\nu,\psi)$ is $mn/2$.

Next consider the transformation
$T_4(\varthetahatvec,\varphihatvec,\varnu)
 = (\thetahatvec,\phihatvec,\nu)$ defined
by~\eqref{varrels}.  The matrix of partial derivatives can be
obtained by substituting \eqref{mnvars} into~\eqref{varrels}.
Without loss of generality, we can suppose that $x_m,y_n=\OO(1)$.
Recall that the Frobenius norm of a matrix is the square root
of the sum of squares of absolute values of the entries.
After multiplying  by $n^{1/2}$ the row indexed by $\nu$ and
dividing by $n^{1/2}$ the column indexed by $\varnu$ (these two
operations together not changing the determinant), the Frobenius
norm of the matrix is $\OO(n^{-1/2})$.  Since the Frobenius norm
bounds the eigenvalues, we can apply
Lemma~\ref{l:detexp} to find that the Jacobian of this
transformation is $1+\OO(n^{-1/2})$.

The transformation $T_4$ changes the region of integration only by
a factor $1+\OO(n^{-1/2})$ in each direction, since the inverse of \eqref{varrels}
has exactly the same form except that the constants $\{ d_j\}, \{d'_k\}$,
while still of magnitude $\OO(n^{-3/2})$, may be different.
Therefore, the image of
region $\X$ lies inside the region
\[
  \Y = \bigl\{\, (\varthetahatvec,\varphihatvec,\varnu) \bigm|
  \abs{\hat\vartheta_j},\abs{\hat\varphi_k}\le 31\kappa~(1\le j\le m, 1\le k\le n),
   \,\abs\varnu\le 31\kappa \, \bigr\}.
\]

We next bound the value of the integrand in $\Y$.   By repeated application
of the inequality $xy\le \tfrac12 x^2+\tfrac12 y^2$, we find that
\[
\dfrac1{12} \sum_{j=1}^m\sum_{k=1}^n A_{jk} (\that_j+\phat_k+\nu)^4
 \le \dfrac{23}{10} \Bigl(\, \sum_{j=1}^m A_{j\b}\hat\vartheta_j^4
    + \sum_{k=1}^n A_{\b k}\hat\varphi_k^4
    + A_{\b\b}\varnu^4\Bigr),
\]
where we have chosen $\tfrac{23}{10}$ as a convenient value
greater than $\tfrac94$ (to cover the small variations in the coefficients)
and less than $\tfrac73$ (to allow us to use Lemma~\ref{ibnd}).
Now define $h(z) = -z^2+\tfrac{23}{10} z^4$.  Then,
for $(\varthetahatvec,\varphihatvec,\varnu)\in\Y$, 
\begin{align}
  \abs{F(\thetavec,\phivec)}
  &\le \exp\Bigl(\, \sum_{j=1}^m A_{j\b|\barXX\,} h(\hat\vartheta_j)
    + \sum_{k=1}^n A_{\b k|\barXX\,} h(\hat\varphi_k) 
    + A_{\b\b|\barXX\,} h(\varnu) \Bigr)
    \notag\\
  &\le \exp\Bigl(\, \sum_{j=1}^{m-1} A_{j\b|\barXX\,} h(\hat\vartheta_j)
    + \sum_{k=1}^{n-1} A_{\b k|\barXX\,} h(\hat\varphi_k) 
    + A_{\b\b|\barXX\,} h(\varnu) \Bigr)
    \label{I199}\\[-0.5ex]
  &= \exp\(A_{\b\b|\barXX\,}h(\varnu)\)
           \prod_{j=1}^{m-1} \exp\(A_{j\b|\barXX\,}h(\hat\vartheta_j)\)
          \prod_{k=1}^{n-1} \exp\(A_{\b k|\barXX\,}h(\hat\varphi_k)\),
     \label{I200}
\end{align}
where the second line holds
because $h(z)\le 0$ for $\abs{z}\le 31\kappa$.

Define
 \begin{align*}
   \W_0 &= \bigl\{ (\varthetahatvec,\varphihatvec,\varnu)\in\Y \bigm|
           \abs{\hat\vartheta_j}\le \tfrac12 n^{-1/2+\eps}
                \quad (1\le j\le m{-}1),\\
           &\kern35mm\abs{\hat\varphi_k}\le \tfrac12 m^{-1/2+\eps}
                    \quad (1\le k\le n{-}1),\\
           &\kern37mm\abs{\varnu}\le \tfrac12 (mn)^{-1/2+2\eps}  
          \bigr\}, \displaybreak[1]\\[0.4ex]
    \W_1&=\Y - \W_0, \\
    \W_2 &= \Bigl\{ \, (\varthetahatvec,\varphihatvec,\varnu)\in\Y \Bigm|\,
                      \Bigl| \, \sum_{j=1}^{m-1} d_j\hat\vartheta_j
                      + \sum_{k=1}^{n-1} d'_k\hat\varphi_k\Bigr|
                                         \le n^{-5/4} \,\Bigr\}.
\end{align*}
Also define similar regions $\W'_0,\W'_1,\W'_2$ by omitting the variables
$\hat{\vartheta}_1,\hat{\varphi}_1$ instead of
$\hat{\vartheta}_m,\hat{\varphi}_n$ starting at~\eqref{I199}.
(Note that without loss of generality we can also assume
that $x_1,\, y_1 = \OO(1)$.)
Using \eqref{varrels}, we see that $T_4$, and the corresponding
transformation that omits $\hat{\vartheta}_1$ and $\hat{\varphi}_1$,
map $\R$ to a superset of
$\W_0\cap\W_2\cap\W'_0\cap\W'_2$.  Therefore,
$\X-\R$ is mapped to a subset of
$\W_1\cup(\W_0-\W_2)\cup\W'_1\cup(\W'_0-\W'_2)$ and it will
suffice to find a tight bound on the integral in each of the four
latter regions.
 
Denoting the right side of~\eqref{I200} by
$F_0(\varthetahatvec,\varphihatvec,\varnu)$,
Lemma~\ref{ibnd} gives
\begin{equation}\label{I0F0}
\int_\Y
    F_0(\varthetahatvec,\varphihatvec,\varnu)
    \, d\varthetahatvec d\varphihatvec d\varnu
    = \exp\(O(m^\eps+n^\eps)\) I_0.
\end{equation}
Also note that
\begin{equation}\label{1Dtail}
 \int_{z_0}^{31\kappa} \exp( c_0 h(z)) = O(1)\exp(c_0 h(z_0))
\end{equation}
for $c_0,z_0>0$ and $z_0=o(1)$, since
$h(z)\le h(z_0)$ for $z_0\le z\le 31\kappa$.
By applying \eqref{1Dtail} to each of the factors of~\eqref{I200} in turn,
\begin{equation}\label{W1bnd}
  \int_{\W_1}  F_0(\varthetahatvec,\varphihatvec,\varnu)
     \, d\varthetahatvec d\varphihatvec d\varnu = 
    O\(e^{-c_6A m^{2\eps} } + e^{-c_6A n^{2\eps} }  \) I_0
\end{equation}
for some $c_6>0$ and so, by \eqref{I0F0} and \eqref{W1bnd},
\[
  \int_{\W_0}  F_0(\varthetahatvec,\varphihatvec,\varnu)
     \, d\varthetahatvec d\varphihatvec d\varnu = 
    \exp\(O(m^\eps+n^\eps)\) I_0.
\]
Applying Lemma~\ref{hoeffding} twice, once to the variables
$d_1\hat\vartheta_1,\ldots,d_{m-1}\hat\vartheta_{m-1},
  d'_1\hat\varphi_1,\ldots,d'_{n-1}\hat\varphi_{n-1}$
and once to their negatives,
using $M=\OO(n^{-2})$, $N=m+n-2$ and $t=n^{-5/4}$,
we find that
\begin{align}
  \int_{\W_0-\W_2} F_0(\varthetahatvec,\varphihatvec,\varnu)
     \, d\varthetahatvec d\varphihatvec d\varnu
     &= O\( e^{-n^{1/4}}\) \int_{\W_0}  F_0(\varthetahatvec,\varphihatvec,\varnu)
     \, d\varthetahatvec d\varphihatvec d\varnu\notag\\
     &= O\( e^{-n^{1/5}}\) I_0.\label{W0W2}
\end{align}

Finally, parallel computations give the same bounds on the integrals
over $\W'_1$ and $\W'_0-\W'_2$.

We have now bounded
$\int \abs{F(\thetavec,\phivec)}$ in regions that together cover the
complement of~$\R$.  Collecting these bounds from
\eqref{R1C}, \eqref{I2mn}, \eqref{W1bnd},
\eqref{W0W2}, and the above-mentioned analogues of
\eqref{W1bnd} and \eqref{W0W2}, we conclude that
\[
 \int_{\R^c} \abs{F(\thetavec,\phivec)}\,d\thetavec d\phivec
=O\(e^{-c_7A m^{2\eps} } + e^{-c_7A n^{2\eps} }  \) I_0
\]
for some $c_7>0$, which Lemma~\ref{boxing} by~\eqref{I0I1}.
\end{proof}

\nicebreak


\begin{thebibliography}{99}

\bibitem{barvinok} A.~Barvinok,
On the number of matrices and a random matrix with prescribed
row and column sums and 0-1 entries,
preprint (2008), available at
\verb+http://arxiv.org/abs/0806.1480+.


\bibitem{Bender} E.\,A.~Bender,
  The asymptotic number of non-negative integer matrices with
  given row and column sums, {\it Discrete Math.}, {\bf 10}
  (1974) 217--223.

\bibitem{Bollobas} B.~Bollob\'as and B.\,D.~McKay,
  The number of matchings in random regular graphs and bipartite
  graphs, {\it J. Combin. Theory Ser.~B},
  {\bf 41} (1986) 80--91.

\bibitem{CGM} E.\,R.~Canfield, C.~Greenhill and B.\,D.~McKay,
  Asymptotic enumeration of dense 0-1 matrices with
  specified line sums, {\it J. Combin. Theory Ser.~A},
  {\bf 115} (2008) 32--66.
  
\bibitem{CM} E.\,R.~Canfield and B.\,D.~McKay,
  Asymptotic enumeration of dense 0-1 matrices with equal row sums
  and equal column sums, {\it Electron. J. Combin.} {\bf 12} (2005), \#R29.

\bibitem{CFM} C.~Cooper, A.~Frieze and M.~Molloy,
  Hamilton cycles in random regular digraphs,
  {\it Combin. Probab. Comput.} {\bf 3} (1994) 39--49. 

\bibitem{GaoTourn} Z.~Gao, B.\,D.~McKay and X.~Wang,
 Asymptotic enumeration of tournaments
 with a given score sequence containing a specified digraph,
 {\it Random Structures and Algorithms}, {\bf 16} (2000) 47--57.
 
\bibitem{GMW} C.~Greenhill, B.\,D.~McKay and X.~Wang,
 Asymptotic enumeration of sparse 0-1 matrices with
irregular row and column sums,
 {\it J. Combin. Theory Ser.~A}, {\bf 113} (2006) 291--324.

\bibitem{Gurvits}
L.~Gurvits,
Van der Waerden/Schrijver-Valiant like conjectures and
stable (aka hyperbolic) homogeneous polynomials: one theorem for
all, \textit{Electron. J. Combin.} \textbf{15} (2008) \#R66
(26 pages).

\bibitem{hoeffding} W.~Hoeffding,
 Probability inequalities for sums of bounded random variables,
 {\it J. Amer. Statist. Assoc.}, {\bf 58} (1963) 13--30.

\bibitem{McKay81} B.\,D. McKay,
 Subgraphs of random graphs with specified degrees,
 {\it Congr. Numer.}, {\bf 33} (1981) 213--223.

\bibitem{Silver} B.\,D.~McKay,
  Asymptotics for 0-1 matrices with prescribed line sums, 
 {\it in } Enumeration and Design, (Academic Press, 1984) 225--238.

\bibitem{RT} B.\,D.~McKay,
 The asymptotic numbers of regular tournaments, eulerian 
 digraphs and eulerian oriented graphs,
 {\it Combinatorica}, {\bf 10} (1990) 367--377.

\bibitem{RANX} B.\,D.~McKay,
 Subgraphs of dense random graphs with specified degrees,
 to appear.

\bibitem{eulerian} B.\,D.~McKay and R.\,W.~Robinson,
Asymptotic enumeration of Eulerian 
 circuits in the complete graph, {\it Combin. Prob. Comput.},
 {\bf 7} (1998) 437--449.

\bibitem{MWtournament} B.\,D.~McKay and X.~Wang,
  Asymptotic enumeration of tournaments with
  a given score sequence, {\it J. Combin. Theory Ser.~A},
  {\bf 73} (1996) 77--90.

\bibitem{MW} B.\,D.~McKay and X.~Wang,
 Asymptotic enumeration of 0-1
 matrices with equal row sums and equal column sums,
 {\it Linear Algebra Appl.}, {\bf 373} (2003) 273--288.

\bibitem{Timashev} A.\,N.~Timash\"ev,
 On permanents of random doubly stochastic matrices and on
 asymptotic estimates for the number of Latin rectangles and
 Latin squares (Russian), {\it Diskret. Mat.}, {\bf 14} (2002)
 65--86; translation in {\it Discrete Math. Appl.}, {\bf 12}
 (2002) 431--452.

\bibitem{Jimmy92} X.~Wang,
Asymptotic enumeration of Eulerian digraphs with multiple edges,
{\it Australas. J. Combin.} {\bf 5} (1992) 293--298.

\bibitem{Jimmy95} X.~Wang,
Asymptotic enumeration of digraphs by excess sequence,
{\it in\/} Graph theory, combinatorics, and algorithms, Vol. 1, 2
(Wiley-Intersci. Publ., Wiley, New York, 1995) 1211--1222.

\bibitem{Jimmy97} X.~Wang,
The asymptotic number of Eulerian oriented graphs with multiple edges,
{\it J. Combin. Math. Combin. Comput.} {\bf 24} (1997) 243--248.

\bibitem{Wormald} N.\,C.~Wormald, Some problems in the enumeration of
  labelled graphs, {\it Ph.\,D.~Thesis}, Department of Mathematics,
  University of Newcastle (1978).

\end{thebibliography}
\end{document}